\documentclass[12pt]{article}
\usepackage{amsmath}
\usepackage{amssymb}
\usepackage{amscd}
 \newcounter{cdef}%
 \newcounter{cthm}%
 \newcounter{cbsp}%
\def\thecdef{\arabic{cdef}}%
\def\thecthm{\arabic{cthm}}%
\def\thecsatz{\arabic{cthm}}%
\def\theclem{\arabic{cthm}}%
\def\thecfolg{\arabic{cthm}}%
\def\thecbsp{\arabic{cbsp}}%
\def\bewname{Proof}%
\def\definname{Definition}%
\def\thmname{Theorem}%
\def\satzname{Proposition}%
\def\lemmaname{Lemma}%
\def\folgname{Corollary}%
\def\bemname{Remark}%
\def\bemsname{Remarks}%
\def\bspname{Example}%
\newcommand{\btxandshort}[1]{and}%
\newcommand{\btxpagesshort}[1]{pp.}%
\newcommand{\Btxinshort}[1]{In}%
\newcommand{\btxphdthesis}[1]{phd-thesis}%
\newcommand{\btxeditorshort}[1]{Ed.}%
\newcommand{\btxeditorsshort}[1]{Eds.}%
\newcommand{\btxvolumeshort}[1]{vol.}%
\newcommand{\btxofseriesshort}[1]{ser.}%

\newenvironment{bew}[1][{}]{\par\addvspace{.5\baselineskip}\par
\textbf{\bewname{#1. }}}{\ \hspace*{\fill}\rule{1ex}{1ex}
\par\addvspace{.5\baselineskip}\par}

\newenvironment{defin}{\par\addvspace{.5\baselineskip}\par
\refstepcounter{cdef} \textbf{\definname~\thecdef.} }{
\par\addvspace{.5\baselineskip}\par}

\newenvironment{thm}{\par\addvspace{.5\baselineskip}\par
\refstepcounter{cthm} \textbf{\thmname~\thecthm.} \em }{\em
\par\addvspace{.5\baselineskip}\par}

\newenvironment{satz}{\par\addvspace{.5\baselineskip}\par
\refstepcounter{cthm} \textbf{\satzname~\thecsatz.} \em }{\em
\par\addvspace{.5\baselineskip}\par}

\newenvironment{lemma}{\par\addvspace{.5\baselineskip}\par
\refstepcounter{cthm} \textbf{\lemmaname~\theclem.} \em }{\em
\par\addvspace{.5\baselineskip}\par}

\newenvironment{folg}{\par\addvspace{.5\baselineskip}\par
\refstepcounter{cthm} \textbf{\folgname~\thecfolg.} \em }{\em
\par\addvspace{.5\baselineskip}\par}

\newenvironment{bem}{\par\addvspace{.5\baselineskip}\par
\textbf{\bemname.} }{\ \hspace*{\fill}\rule{1ex}{1ex}
\par\addvspace{.5\baselineskip}\par}

\newenvironment{bems}{\par\addvspace{.5\baselineskip}\par
\textbf{\bemsname.} }{\ \hspace*{\fill}\rule{1ex}{1ex}
\par\addvspace{.5\baselineskip}\par}

\newenvironment{bsp}{\par\addvspace{.5\baselineskip}\par
\refstepcounter{cbsp}
\textbf{\bspname~\thecbsp.} }{\ \hspace*{\fill}\rule{1ex}{1ex}
\par\addvspace{.5\baselineskip}\par}

\newcommand{\linv}[1]{(#1)_{\mathrm{l}}}

\newcommand{\dif}{\mathrm{d}}

\newcommand{\id}{\mathrm{id}}

\newcommand{\Lin}{\mathrm{Lin}}

\newcommand{\mc}{\mathcal}
\newcommand{\SLq}[1]{\mathrm{SL}_q(#1)}
\newcommand{\SUq}[1]{\mathrm{SU}_q(#1)}
\newcommand{\OSLq}[1]{\mathcal{O}(\mathrm{SL}_q(#1))}
\newcommand{\dOSLq}[1]{\mathcal{O}(\mathrm{SL}_q(#1))^\circ }
\newcommand{\Ga}{\varGamma }
\newcommand{\uGa}{{}_\mathrm{u}\Ga }
\newcommand{\A}{\mathcal{A}}
\newcommand{\Anull}{\A ^\circ }
\newcommand{\copr}{\varDelta }
\newcommand{\coun}{\varepsilon }
\newcommand{\antip}{S}
\newcommand{\lcoa}{\copr _{\scriptscriptstyle \mathrm{L}}}
\newcommand{\ot}{\otimes }
\newcommand{\otA}{\ot _{\scriptscriptstyle \A }}
\newcommand{\w}{\omega }
\newcommand{\wH}{\w _{\scriptscriptstyle H}}
\newcommand{\wX}{\w _{\scriptscriptstyle X}}
\newcommand{\wY}{\w _{\scriptscriptstyle Y}}
\newcommand{\cHX}{c _{\scriptscriptstyle HX}}
\newcommand{\cHY}{c _{\scriptscriptstyle HY}}
\newcommand{\cXY}{c _{\scriptscriptstyle XY}}
\newcommand{\x}{X}
\newcommand{\X}{\mathcal{X}}
\newcommand{\comp}{\mathbb{C}}
\newcommand{\compx}{\mathbb{C}^\times }
\renewcommand{\linv}[1]{{#1}_{\scriptscriptstyle \mathrm{L}}}
\newcommand{\paar}[2]{\langle #1,#2 \rangle }
\newcommand{\E}[1]{E^{(#1)}}
\newcommand{\F}[1]{F^{(#1)}}
\newcommand{\G}[1]{G^{(#1)}}
\newcommand{\ii}{\mathrm{i}}

\title{Classification of Left-Covariant Differential Calculi on the Quantum
Group $\SLq 2$\thanks{This paper was supported by the
Deutsche Forschungsgemeinschaft}}
\author{Istv\'an Heckenberger\thanks{e-mail:
heckenbe@mathematik.uni-leipzig.de}\\
{\small Universit\"at Leipzig, Augustusplatz 10-11,}\\
{\small 04109 Leipzig, Germany}}

\begin{document}
\maketitle
\begin{abstract}
For transcendental values of $q$ the quantum tangent spaces of all
left-covariant first order differential calculi of dimension less than
four on the quantum group $\SLq 2$ are given.
All such differential calculi $\Ga $ are determined and investigated
for which the left-invariant differential one-forms $\omega (u^1_2)$,
$\omega (u^2_1)$ and $\omega (u^1_1-u^2_2)$ generate $\Ga $ as a bimodule
and the universal
higher order differential calculus has the same dimension as in the classical
case. Important properties (cohomology spaces, $*$-structures, braidings,
generalized Lie brackets) of these calculi are examined as well.
\end{abstract}

\noindent
\textbf{MSC (1991):} 17B37, 46L87, 81R50\\
\textbf{Keywords:} quantum groups, noncommutative differential calculus,
quantum tangent space

\section{Introduction}
\label{s-intro}

The theory of bicovariant differential calculus over Hopf algebras
is one of the commonly used and best understood theories related to
quantum groups. Its origin was a paper of S.\,L.~Woronowicz \cite{a-Woro2}
where also left-covariant differential calculi were considered.

It is well known that bicovariant differential calculi on quantum groups often
(but not always, see \cite{a-Karim1},\cite{a-CornJac1}) have the
unpleasant property (besides non-uniqueness) that their dimensions do not
coincide with the dimensions of the canonical differential calculi of the
corresponding Lie groups. There were made some attempts using a generalized
adjoint action in order to circumvent this defect
\cite{a-DelHueff1},\cite{a-Sudbery1}.
Another way is to look for left-covariant differential
calculi on the quantum group. Then the corresponding quantum tangent spaces
are not invariant under the adjoint action in general.
The first such example (the legendary 3D-calculus)
was developed by S.\,L.~Woronowicz \cite{a-Woro3} for the quantum group
$\SUq 2$.
Among others it was shown therein that the cohomology spaces of the differential
complex are the same as in the classical situation.
Further examples of such kind were given by K.~Schm\"udgen and A.~Sch\"uler
\cite{a-SchSch3,inp-SchSch2}.
The paper \cite{a-SchSch3} contains also a first classification of
left-covariant
differential calculi on the quantum group $\SLq 2$ (under very restrictive
conditions). A method for the construction of left-covariant differential
calculi on a quantum linear group was initiated by K.~Schm\"udgen
\cite{a-Schmue1}.
Since the pioneering work of Woronowicz the general theory of differential
calculi on quantum groups was refined and developed further.
An extensive overview can be found in Chapter 14 of the monograph
\cite{b-KS}.

In contrast to the classical situation there is no distinguished
differential calculus on a quantum group and the non-commutative geometry of
a quantum group depends on the differential calculus in general.
Thus it seems to be natural to ask how many such calculi
(satisfying other additional natural conditions) do really exist.
This problem is studied in the present paper. Our aim is to give a
step-by-step classification of left-covariant
differential calculi $(\Ga ,\dif )$ on the quantum group $\SLq 2$ under
rather general assumptions. In order
to motivate the addition of further assumptions we investigate the outcome of
the classification after each step. Although the methods we use can be easily
described, the computations are very boring. The
computer algebra program FELIX \cite{inp-ApelKlaus} by J.~Apel and U.~Klaus
was very helpful to carry out long computations.
The main result of the present paper is Theorem \ref{t-main}. Suppose that $q$
is a nonzero complex number and not a root of unity. Then the assertion of
Theorem \ref{t-main} states that there are exactly
11 (single) left-covariant first order differential calculi over the
Hopf algebra $\OSLq 2$ having the following properties:
The quantum tangent space
is a subspace of the algebra $\mc{U}$ (see Section \ref{sec-qgslq2}),
the one-forms $\omega (u^1_2)$, $\omega (u^2_1)$ and $\omega (u^1_1-u^2_2)$
form a basis of the (necessarily 3-dimensional)
$\OSLq 2$-bimodule $\Ga $, the dimension of the space of left-invariant
differential 2-forms in the universal higher order differential calculus is
at least 3 and the first order calculus is invariant with respect to all
Hopf algebra automorphisms of $\OSLq 2$. The list of these calculi is given
in Corollary \ref{f-Hopfinvcal}. Using the method of Woronowicz
\cite{a-Woro3} it is proved in Theorem \ref{t-cohom}
that the dimensions of the cohomology spaces of these 11
differential complexes are the same as in the classical situation.

This paper is organized as follows. In Section~\ref{sec-lcdc} we recall some
basic notions and facts about the general theory of left-covariant differential
calculus on quantum groups. If not stated otherwise we follow the definitions
and notations of Woronowicz \cite{a-Woro2} and of the monograph
\cite{b-KS}. In Section~\ref{sec-qgslq2} the structure of the dual
Hopf algebra $\mc{U}$ of $\OSLq 2$ is described. In Section~\ref{sec-urcoid}
we determine all 4-dimensional unital right coideals of $\mc{U}$.
In Section~\ref{sec-restr} further restrictions on the calculus are added.
In Section~\ref{sec-struc} we investigate additional structures
such as $*$-structures and braidings. In Section~\ref{sec-deRham}
the cohomology spaces of the most important differential complexes
are studied.
Section~\ref{sec-detail} contains the main theorem
(Theorem~\ref{t-main}) of the present paper. The outcoming calculi
are then studied in detail.

Throughout Sweedler's notation for coproducts and coactions and
Einsteins convention of summing over repeated indices are used.
The symbols $\ot $ and $\otA $ denote tensor products over the complex numbers
and over an algebra $\A $, respectively. All algebras are complex and unital.

\section{Left-covariant differential calculi on quantum groups}
\label{sec-lcdc}

First let us recall some facts of the general theory (see \cite{a-Woro2},
\cite{b-KS}). Let $\A $ be a Hopf algebra with coproduct $\copr $,
counit $\coun $, and invertible antipode $\antip $.
An $\A $-bimodule $\Ga $ is called \textit{first order differential
calculus} (FODC for short) over $\A $, if there is a linear mapping
$\dif :\A \to \Ga $ such that
\begin{itemize}
\item
$\dif $ satisfies the Leibniz rule: $\dif (ab)=(\dif a)b+a\dif b$ for any
$a,b\in \A $,
\item
$\Ga =\Lin \{a\dif b\,|\,a,b\in \A \}$.
\end{itemize}
An FODC $\Ga $ is called \textit{left-covariant} if there is a linear mapping
$\lcoa :\Ga \to \A \ot \Ga $ such that
$\lcoa (a(\dif b)c)=\copr (a)\cdot (\id \ot \dif )\copr (b)\cdot \copr (c)$,
where $(a\ot b)\cdot (c\ot \rho )=ac\ot b\rho $ and
$(a\ot \rho )\cdot (b\ot c)=ab\ot \rho c$ for any $a,b,c\in \A $ and
$\rho \in \Ga $.
Elements $\rho \in \Ga $ for which $\lcoa (\rho )=1\ot \rho $ are called
\textit{left-invariant}.
Because of Theorem 2.1 in \cite{a-Woro2} any left-covariant
$\A $-bimodule is a free left module and any basis of the vector space
$\linv{\Ga }$ of left-invariant 1-forms is a free basis of the
left (right) $\A $-module $\Ga $. The dimension of $\linv{\Ga }$ is called
the \textit{dimension} of the FODC $\Ga $. In this paper we are dealing
only with finite dimensional FODC.

Suppose that $\Ga $ is an $n$-dimensional first order differential calculus
over $\A $. Let us fix a basis $\{\w _i\,|\,i=1,\ldots ,n\}$ of
$\linv{\Ga }$. Then there are functionals $\x _i$, $i=1,\ldots ,n$
in the \textit{dual Hopf algebra} $\Anull $
such that the differential $\dif $ can be written in the form
\begin{align}\label{eq-difa}
\dif a=\sum _{i=1}^na_{(1)}\x _i(a_{(2)})\w _i, \quad a\in \A .
\end{align}
Recall that $\Anull $ is the set of all linear functionals $f$ on $\A $
for which there exist functionals $f_1,\ldots,f_N, g_1,\ldots, g_N$ on $\A$
such that $f(ab)=\sum _{i=1}^Nf_i(a)g_i(b)$ for all $a,b\in \A $.

The vector space $\X _\Ga :=\Lin \{\x _i\,|\,i=1,\ldots ,n\}$ is called the
\textit{quantum tangent space} of the left-covariant FODC $\Ga $.
We define a mapping $\omega :\A \to \linv{\Ga }$ by
$\omega (a)=S(a_{(1)})\dif a_{(2)}$. Then by (\ref{eq-difa}) the equation
\begin{align}\label{eq-omegaa}
\omega (a)=\sum _{i=1}^n\x _i(a)\w _i, \quad a\in \A
\end{align}
holds. Since $\dif 1=0$, we have $\omega (1)=0$ and
$\x _i(1)=0$ for any $i$. The following lemma is the starting point for the
first part of our classification.

\begin{lemma}\label{l-lcdc}
If $\X $ is the quantum tangent space of a FODC $\Ga $, then
$\bar{\X }=\X \oplus \comp \coun $ is a unital right coideal of $\Anull $
(i.\,e.\ $\copr (\bar{\X }) \subset \bar{\X } \ot \Anull $).\\
Conversely, any unital right coideal $\bar{\X }$ of $\Anull $ determines
a unique FODC with quantum tangent space
$\X ^+:=\{X \in \bar{\X }\,|\,X(1)=0\}$.
\end{lemma}

\begin{bew}
See \cite{b-KS} and \cite{a-HeckSchm2}.
\end{bew}

In particular, the coproduct of elements of the quantum tangent space takes
the form
\begin{align}
\copr \x _i&=1\ot \x _i+\x _j\ot f^j_i,
\end{align}
where the functionals $f^j_i\in \Anull $ describe the bimodule structure
of $\Ga $:
\begin{align}\label{eq-module}
\w _ia&=a_{(1)}f^i_j(a_{(2)})\w _j \qquad \text{for any $a\in \A $.}
\end{align}

Left-covariant first order differential calculi $\Ga $ over $\A $ are also
characterized by the right ideal
\begin{gather}\label{eq-rightideal}
\mc{R}_\Ga :=\{a\in \ker \coun \subset \A \,|\,\omega (a)=0 \}
=\{a \in \A \,|\,X(a)=0\quad \forall X\in \bar{\X }_\Ga \}
\end{gather}
of $\A $. Two FODC $(\Ga _1,\dif _1)$ and $(\Ga _2,\dif _2)$ over $\A $
are called \textit{isomorphic}, if $\mc{R}_{\Ga _1}=\mc{R}_{\Ga _2}$
or equivalently if $\X _{\Ga _1}=\X _{\Ga _2}$.

Let $\Ga ^{\ot k}$ denote the $k$-fold tensor product
$\Ga \otA \cdots \otA \Ga $ of the $\A $-bimodule $\Ga $,
$\Ga ^{\ot 0}:=\A $, $\Ga ^{\ot 1}:=\Ga $ and
$\Ga ^\ot :=\bigoplus _{k=0}^\infty \Ga ^{\ot k}$. Then $\Ga ^\ot $ becomes
an algebra with multiplication $\otA $. Let $\mc{S}$ be a graded two-sided
ideal in $\Ga ^\ot $, $\mc{S}\subset \bigoplus _{k=2}^\infty \Ga ^{\ot k}$,
$\mc{S}=\bigoplus _{k=2}^\infty \mc{S}\cap \Ga ^{\ot k}$.
Then the $\A $-bimodule $\Ga ^\wedge :=\Ga ^\ot /\mc{S}$ as well as
$\Ga ^{\wedge k}:=\Ga ^{\ot k}/(\mc{S}\cap \Ga ^{\ot k})$ are well defined.
The bimodule $\Ga ^\wedge $ is called a \textit{differential calculus} over
the Hopf algebra $\A $ if there is a linear mapping
$\dif :\Ga ^\wedge \to \Ga ^\wedge $ of
grade one (i.\,e.\ $\dif :\Ga ^{\wedge k}\to \Ga ^{\wedge k+1}$) such that
\begin{itemize}
\item
$\dif $ satisfies the graded Leibniz rule $\dif (\rho \wedge \rho ')=
\dif \rho \wedge \rho '+(-1)^m\rho \wedge \dif \rho '$ for
$\rho \in \Ga ^{\wedge m},\rho '\in \Ga ^\wedge $,
\item
$\dif ^2=0$,
\item
$\Ga = \Lin \{a\dif b\,|\,a,b\in \A \}$.
\end{itemize}
If $\Ga $ is left-covariant, then $\Ga ^\ot $ is also left-covariant with
$\lcoa (\rho \otA \rho ')=\rho _{(-1)}\rho '_{(-1)}\ot
\rho _{(0)}\otA \rho '_{(0)}$. Suppose that $\mc{S}$ is an invariant subspace
of the left coaction, i.\,e.\ $\lcoa (\mc{S})\subset \A \ot \mc{S}$.
Then $\Ga ^\wedge $ inherits the left coaction of $\Ga ^\ot $ and
$\Ga ^\wedge $ is called a \textit{left-covariant differential calculus}
over $\A $.

Suppose that $\Ga ^\wedge $ is a left-covariant differential calculus
over $\A $. Then the Maurer-Cartan formula
\begin{align}
\dif \omega (a)=-\omega (a_{(1)})\wedge \omega (a_{(2)}),\quad a\in \A
\end{align}
is always fulfilled. Moreover, for any given left-covariant FODC $\Ga $
over $\A $ there exists a universal differential calculus
$\uGa ^\wedge $. This means that any left-covariant differential calculus
$\tilde{\Ga }^\wedge $ over $\A $ with $\tilde{\Ga }^{\wedge 1}=\Ga $
is isomorphic to $\uGa ^\wedge /\tilde{\mc{S}}$, where $\tilde{\mc{S}}$
is a two-sided ideal in $\uGa ^\wedge $. The differential calculus
$\uGa ^\wedge $ can be given by the two-sided ideal $\mc{S}$ generated
by the elements of the vector space
\begin{gather}\label{eq-symmform}
\linv{\mc{S}^2}:=
\Lin \{\omega (a_{(1)})\otA \omega (a_{(2)})\,|\,a\in \mc{R} \}.
\end{gather}

\begin{lemma}\label{l-dimGaw}
The following equation holds for any left-covariant differential calculus
$\Ga $ over $\A $ with quantum tangent space $\X $:
\begin{align}
\dim \linv{(\uGa ^{\wedge 2})}=
\dim \{T\in \bar{\X }\ot \bar{\X }\,|\,\mathrm{m}T=0\}-\dim \X,
\end{align}
where $\mathrm{m}$ denotes the multiplication map
$\mathrm{m}:\bar{\X }\ot \bar{\X }\subset \Anull \ot \Anull \to \Anull$.
\end{lemma}

\begin{bew}
It was proved in \cite{a-Schueler2} that in the present situation the formula
\begin{align}
\dim \linv{(\uGa ^{\wedge 2})}
=\dim \{T\in \X \ot \X \,|\,\mathrm{m}T\in \X \}
\end{align}
is valid. Let $V,V'$, and $W$ denote the vector spaces
\[V:=\{T\in \X \ot \X \,|\,\mathrm{m}T\in \X \},\quad
V':=\{T\in \bar{\X }\ot \X \,|\,\mathrm{m}T=0\},\quad \text{and}
\]
\[W:=\{T\in \bar{\X }\ot \bar{\X }\,|\,\mathrm{m}T=0\},
\]
respectively.
Then the mapping $\varphi :V \to V'$, $\varphi (T)=T-1\ot \mathrm{m}T$,
and its inverse $\psi :V'\to V$, $\psi (T)=T -(\coun \ot \id)T$
give an isomorphism between $V$ and $V'$. Obviously we have
$\bar{\X }\ot \bar{\X }=(\bar{\X }\ot \X )\oplus (\bar{\X }\ot \comp \cdot 1)$
and hence
\[\dim W=\dim V'+\dim W',
\]
where $W':=\{T'\in \bar{\X }\,|\,\exists T\in \bar{\X }\ot \X ,
\mathrm{m}(T'\ot 1-T)=0\}$.
But the vector space $W'$ is isomorphic to $\X $. Indeed, $\X \subset W'$
since
$\mathrm{m}(X\ot 1-1\ot X)=0$ but $1\not \in W'$. The latter follows from the
fact that $\coun (\mathrm{m}(1\ot 1))=1$ and
$\coun (\mathrm{m}(T))=0$ for any $T\in \bar{\X }\ot \X $.
\end{bew}

\section{The Hopf dual of {\boldmath $\OSLq 2$}}
\label{sec-qgslq2}

In what follows we assume that $q$ is a transcendental complex number.
The structure of the coordinate Hopf algebra $\OSLq 2$
(with generators $u^i_j$, $i,j=1,2$) of the quantum group $\SLq 2$ is well
known. We now restate the description of the Hopf dual $\mc{U}=\dOSLq 2$
obtained in the monograph \cite{b-Joseph}. Let $\mc{U}$ denote the unital
algebra generated by the elements $E,F,G$ and $f_\mu $
($\mu \in \compx =\comp \setminus \{0\}$) and by the relations
\begin{equation}
\begin{gathered}
f_\mu f_\nu =f_{\mu \cdot \nu },\quad
f_\mu E=\mu ^2 Ef_\mu,\quad
f_\mu F=\mu ^{-2} Ff_\mu,\quad
f_\mu G=Gf_\mu ,\\
GE=E(G+2),\quad GF=F(G-2),\quad EF-FE=\frac{f_q-f_{q^{-1}}}{q-q^{-1}}.
\end{gathered}
\end{equation}
The element $f_1$ is the unit in the algebra $\mc{U}$. The element $f_{-1}$ is
also denoted by $\coun _-$.

Let us fix one square root $q^{1/2}$ of $q$ and define $K=f_{q^{1/2}}$.
Then there is a Hopf algebra structure on $\mc{U}$ such that
\begin{equation}
\begin{aligned}
\copr (E)&=E\ot K+K^{-1}\ot E,\quad &\coun (E)&=0,\quad & \antip (E)&=-qE,\\
\copr (F)&=F\ot K+K^{-1}\ot F,\quad &\coun (F)&=0, & \antip (F)&=-q^{-1}F,\\
\copr (G)&=1\ot G+G\ot 1, &\coun (G)&=0, & \antip (G)&=-G,\\
\copr (f_\mu )&=f_\mu \ot f_\mu , &\coun (f_\mu )&=1, &
\antip (f_\mu )&=f_{\mu ^{-1}}.
\end{aligned}
\end{equation}
To make calculations easier we use the notation
$\F k:=F^kK^{-k}/[k]!$, $\E k:=K^{-k}E^k/[k]!$
and $\G k=G^k/k!$ for any $k\in \mathbb{N}_0$, where
$[k]!=[k][k-1]\cdots [1]$, $[0]!=1$ and $[k]=(q^k-q^{-k})/(q-q^{-1})$.
The coproducts of these elements are
\begin{equation}\label{eq-scopr}
\begin{aligned}
\copr (\F k)& =\sum _{r=0}^k\F{k-r}K^{-2r}\ot \F r,\\
\copr (\E k)& =\sum _{s=0}^kK^{-2s}\E{k-s}\ot \E s,\\
\copr (\G k)& =\sum _{t=0}^k\G{k-t}\ot \G t.
\end{aligned}
\end{equation}

The dual pairing $\paar{\cdot }{\cdot }$ of the Hopf algebras
$\mc{U}=\dOSLq 2$ and $\OSLq 2$ is given by the matrices
$\paar{\cdot }{u^i_j}$, where
\begin{gather}\label{eq-pairing}
E=\begin{pmatrix}0&0\\1&0\end{pmatrix},\quad
F=\begin{pmatrix}0&1\\0&0\end{pmatrix},\quad
G=\begin{pmatrix}-1&0\\0&1\end{pmatrix},\quad
f_\mu =\begin{pmatrix}\mu ^{-1}&0\\0&\mu \end{pmatrix}.
\end{gather}

For the algebra $\mc{U}$ a PBW-like theorem holds:
the elements of the set $\{\F if_\mu \E j\G k\,|\,i,j,k\in \mathbb{N}_0,
\mu \in \compx \}$ form a vector space basis of the algebra $\mc{U}$.

\begin{bem}
That the Hopf algebra $\mc{U}$ defined above is the full Hopf dual of the
Hopf algebra $\OSLq 2$ is proved in the monograph \cite{b-Joseph}
for transcendental $q$. The assumption that $q$ is transcendental is only
necessary in order to apply this theorem. If this fact holds already for
$q$ not a root of unity, then all results of the present paper remain valid
under this assumption. In what follows, we classify the corresponding
left-covariant differential calculi whose quantum tangent space are
contained in $\mc{U}$. All these considerations hold if $q\not= 0$ and
$q$ is not a root of unity.
\end{bem}

\section{Unital right coideals of {\boldmath $\mc{U}$}}
\label{sec-urcoid}

By Lemma \ref{l-lcdc} left-covariant first order differential calculi
over the Hopf algebra $\A :=\OSLq 2$ and unital right coideals
of $\Anull $ are in one-to-one correspondence.
Now we determine all unital right coideals of $\mc{U}$ of dimension
$\leq 4$. In this way we decribe
all left-covariant first order differential calculi over $\A $ of dimension
less than four.

For $X=\F if_\mu \E j\G k$, $i,j,k\in \mathbb{N}_0$, $\mu \in \compx $ we set
\begin{align*}
\partial _1(X):=i,\quad \partial _2(X):=\mu ,\quad
\partial _3(X):=j,\quad \partial _4(X):=k,\quad \partial _{13}(X):=i+j.
\end{align*}
For a finite linear combination $X=\sum_{i=1}^na_iX_i$ of such elements we
define $\partial _m(X):=\mathrm{max}\{\partial _m(X_i)\,|\,i=1,\ldots,n\}$
for $m\in \{1,3,4,13\}$.
Then from the PBW-theorem we conclude that
$\mc{U}=\bigoplus _{\mu \in \compx }[\mc{U}]_\mu $, where
\[
[\mc{U}]_\mu :=\Lin \{\F if_\mu \E j\G k\,|\,i,j,k\in \mathbb{N}_0\}.
\]

\begin{satz}\label{s-basissatz}
Let $\X $ be a right coideal of $\mc{U}$ and let $X\in \X $.
Then there exist complex numbers $\alpha _{i,\nu ,j,k}$
($\nu \in \compx $, $i,j,k\in \mathbb{N}_0$) such that
$X=\sum _{i,\nu ,j,k}\alpha _{i,\nu ,j,k}\F if_\nu \E j\G k$.
Let us fix $\mu \in \compx $ and $r,s,t\in \mathbb{N}_0$.
Consider the element
\begin{align}
X_{r,\mu ,s,t}&:=\sum _{i,j,k}\alpha _{i,\mu ,j,k}
\F{i-r}f_{q^{-r-s}\mu }\E{j-s}\G{k-t}
\end{align}
of $\mc{U}$, where the sum is running over all $i,j,k\in \mathbb{N}_0$
with $i\geq r$, $j\geq s$ and $k\geq t$.
Then the vector space $\Lin \{X_{r\mu st}\,|\,r,s,t\in \mathbb{N}_0,
\mu \in \compx \}$ is the smallest (with respect to inclusion)
right coideal of $\mc{U}$ containing $X$.
\end{satz}

\begin{bew}
Using (\ref{eq-scopr}) we compute the coproduct of $X$ and obtain
\begin{align*}
&\copr \left(\sum _{i\nu jk}\alpha _{i\nu jk}\F if_\nu \E j\G k\right)=\\
&\qquad =
\sum _{i\nu jk}\alpha _{i\nu jk}\sum _{rst}
\F{i-r}K^{-2r}f_\nu K^{-2s}\E{j-s}\G{k-t}\ot \F rf_\nu \E s\G t\\
&\qquad =
\sum _{r\nu st}\sum _{ijk}\alpha _{i\nu jk}
\F{i-r}f_{q^{-r-s}\nu }\E{j-s}\G{k-t}\ot \F rf_\nu \E s\G t\\
&\qquad =
\sum _{r\nu st}X_{r\nu st}\ot \F rf_\nu \E s\G t.
\end{align*}
Because of the PBW-theorem the elements of the right hand side of the tensor
product are linearly independent. Hence the elements $X_{r\nu st}$ belong to
any right coideal containing $X$ (i.\,e.\ they belong to $\X $ too).
Observe that $X=\sum _{\mu \in \compx }X_{0\mu 00}$.
Taking in account that $(\copr \ot \id )\copr (X)=(\id \ot \copr )\copr (X)$
it follows that the vector space
$\Lin \{X_{r\mu st}\,|\,r,s,t\in \mathbb{N}_0, \mu \in \compx \}$
is a right coideal of $\mc{U}$.
\end{bew}

Since $X=\sum _{\mu \in \compx }X_{0\mu 00}$ we obtain the following.

\begin{folg}\label{f-homrc}
Any right coideal $\X $ of $\mc{U}$ is isomorphic to the direct sum
$\bigoplus _{\mu \in \compx }[\mc{X}]_\mu $ of its homogeneous
components $[\X ]_\mu :=\X \cap [\mc{U}]_\mu $.
\end{folg}

What can be said about the dimension of a right coideal $\X $ of $\mc{U}$?
Let $X$ be a nonzero element of $\X $. By Corollary \ref{f-homrc} we may
assume without loss of generality that there is a $\mu \in \compx $ such that
$X\in [\mc{X}]_\mu $. Due to the PBW-theorem there is for any
$t\in \mathbb{N}_0$, $t\leq \partial _4(X)$ a unique $X_t\in \mc{U}$ such that
$X=\sum _{t=0}^{\partial _4(X)}X_t\G t$ and $X_{\partial _4(X)}\not=0$. Let now
$p\in \mathbb{N}_0$, $p\leq \partial _{13}(X)$. We define
\begin{align*}
t_{13}(p)=\mathrm{max}\{m\in \mathbb{Z}\,|\,X_m\not= 0,
\partial _{13}(X_m)\geq p\}.
\end{align*}
Then the coefficients of $X_{r\mu st}$ in Proposition \ref{s-basissatz}
indicate that for any
$t\in \mathbb{N}_0$, $t\leq t_{13}(p)$ there is at least one
number $r(t,p)\in \mathbb{N}_0$, $r(t,p)\leq p$ such that
$X_{r(t,p),\mu ,p-r(t,p),t}\not=0$. Since
$\partial _2(X_{r(t,p),\mu ,p-r(t,p),t})=\mu q^{-p}$, $q$ is not a root of
unity and $\partial _4(X_{r(t,p),\mu ,p-r(t,p),t})=t_{13}(p)-t$
we conclude that the elements
$X_{r(t,p),\mu ,p-r(t,p),t}$, $0\leq p\leq \partial _{13}(X)$,
$0\leq t\leq t_{13}(p)$ are linearly independent.
Hence the number
\begin{align}\label{eq-lbdim}
\sum _{p=0}^{\partial _{13}(X)}(t_{13}(p)+1)
=\sum _{t=0}^{\partial _4(X)}(s_{13}(t)+1)
\end{align}
with $s_{13}(t)=\mathrm{max}\{\partial _{13}(X_m)\,|\,m\geq t, X_m\not= 0\}$
for $t\in \mathbb{N}_0$, $t\leq \partial _4(X)$,
is a lower bound for the dimension of the right coideal $\X $.

\begin{bsp}
Suppose that $X=K^6\G 5+(\F 1 K^6-K^6)\G 4+(\F 2K^6\E 1+K^6\E 3)\G 2+\F 3K^6$
is an element of a right coideal $\X $ of $\mc{U}$. We have
$X\in [\mc{U}]_{q^3}$, $\partial _4(X)=5$ and $\partial _{13}(X)=3$.
The coefficients $X_k$ of $\G k$ in $X$ are $X_0=\F 3K^6$,
$X_2=\F 2K^6\E 1+K^6\E 3$,
$X_4=\F 1K^6-K^6$, $X_5=K^6$ and $X_t=0$ otherwise. The values of the
functions $t_{13}(n)$ ($s_{13}(n)$) are $5,4,2,2$ for $n=0,1,2,3$
($3,3,3,1,1,0$ for $n=0,1,2,3,4,5$).
Hence the dimension of each right coideal $\mc{X}$ containing $X$ is at least
17.
\end{bsp}

In the remaining part of this section we determine all unital right coideals
of $\mc{U}$ of dimension $\leq 4$.

\subsection{{\boldmath $\dim \X \leq 2$}}

The following list contains all possibilities (if not otherwise
stated, parameters are arbitrary complex numbers):
\begin{itemize}
\item
$\X ^1_1=\comp \cdot 1$.
\item
$\X ^2_1=\Lin \{G,1\}$,
\item
$\X ^2_2=\Lin \{\F 1K^2+\alpha K^2\E 1+\beta K^2,1\}$,
\item
$\X ^2_3=\Lin \{K^2\E 1+\alpha K^2,1\}$,
\item
$\X ^2_4=\Lin \{f_\mu ,1\}$, $\mu \in \compx $, $\mu \not= 1$.
\end{itemize}

Obviously $\dim \X =1$ and $1\in \X $ imply that the only 1-dimensional
unital right coideal of $\mc{U}$ is $\X ^1_1$.

Suppose that $\dim \X =2$.
By (\ref{eq-lbdim}) we have $\partial _4(X)<\dim \X =2$ for any $X\in \X $.
If there is an $X\in [\X ]_\mu $, $\mu \in \compx $ with $\partial _4(X)=1$
then we must have $s_{13}(0)=s_{13}(1)=0$ by (\ref{eq-lbdim}) and we obtain
$X=\alpha f_\mu G+\beta f_\mu $, $\alpha \not= 0$.
The only nonzero elements $X_{r\nu st}$ are $X$ and
$X_{0\mu 01}=\alpha f_\mu $.
Since $\dim \X =2$ and $1\in \X $, $1$ must be one of those two elements.
Therefore, $f_\mu =1$, i.\,e.\ $\mu =1$.
Hence $\X $ is isomorphic to $\X ^2_1$.

If $\partial _4(X)=0$ for any $X\in \X $ then by (\ref{eq-lbdim}) $s_{13}(0)$
can take the
values $1$ and $0$. In the first case we have $1=s_{13}(0)=\partial _{13}(X)$
from which $X=\alpha \F 1f_\mu +\beta \E 1f_\mu +\gamma f_\mu $
($\alpha \not=0$ or $\beta \not=0$) follows. Then
$X_{1\mu 00}=\alpha K^{-2}f_\mu $, $X_{0\mu 01}=\beta K^{-2}f_\mu $
and therefore $K^{-2}f_\mu \in \X $. Since $1\in \X $, we must have
$f_\mu =K^2$.
This gives the coideals $\X ^2_2$ and $\X ^2_3$.
In the remaining case we have $\partial _{13}(X)=0$ and we get $\X ^2_4$.

\subsection{{\boldmath $\dim \X =3$}}

By (\ref{eq-lbdim}), $\partial _4(X)\leq 2$ for any $X\in \X $.
If $\partial _4(X)=2$ then $s_{13}(i)=0$ for any $i=0,1,2$.
If $\partial _4(X)=1$ then $s_{13}(1)=0$ and $s_{13}(0)$ can take the values
$0$ and $1$.
Finally, if $\partial _4(X)=0$ then $s_{13}(0)$ can be $2,1$ or $0$.
We obtain the following list of unital right coideals:

\begin{itemize}
\item
$\X ^3_1=\Lin \{G^2,G,1\}$,
\item
$\X ^3_2=\Lin \{K^2G+\alpha \F 1K^2+\beta K^2\E 1,K^2,1\}$,
$|\alpha |+|\beta |\not= 0$,
\item
$\X ^3_3=\Lin \{G+\alpha \F 1+\beta \E 1,K^{-2},1\}$,
$|\alpha |+|\beta |\not= 0$,
\item
$\X ^3_4=\Lin \{f_\mu G,f_\mu ,1\}$, $\mu \in \compx $, $\mu \not=1$,
\item
$\X ^3_5=\Lin \{G,\F 1K^2+\alpha K^2\E 1+\beta K^2,1\}$,
\item
$\X ^3_6=\Lin \{G,K^2\E 1+\alpha K^2,1\}$,
\item
$\X ^3_7=\Lin \{G,f_\mu ,1\}$, $\mu \in \compx $, $\mu \not=1$,
\item
$\X ^3_8=\Lin \{\F 2K^4+\alpha \F 1K^4\E 1+\alpha ^2K^4\E 2
+\beta \F 1K^4+\alpha \beta K^4\E 1+ \gamma K^4,
\F 1K^2+\alpha K^2\E 1+\beta K^2,1\}$,
\item
$\X ^3_9=\Lin \{K^4\E 2+\alpha K^4\E 1+\beta K^4,K^2\E 1+\alpha K^2,
1\}$,
\item
$\X ^3_{10}=\Lin \{\F 1f_\mu +\alpha f_\mu \E 1+\beta f_\mu ,
f_{\mu q^{-1}},1\}$, $\mu \not= q$,
\item
$\X ^3_{11}=\Lin \{f_\mu \E 1+\alpha f_\mu ,f_{\mu q^{-1}},1\}$,
$\mu \not= q$,
\item
$\X ^3_{12}=\Lin \{\F 1K^2+\alpha K^2, K^2\E 1+\beta K^2,1\}$,
\item
$\X ^3_{13}=\Lin \{\F 1K^2+\alpha K^2\E 1+\beta K^2, f_\mu ,1\}$,
$\mu \not= 1$,
\item
$\X ^3_{14}=\Lin \{K^2\E 1+\alpha K^2, f_\mu ,1\}$, $\mu \not= 1$,
\item
$\X ^3_{15}=\Lin \{f_\mu ,f_\nu ,1\}$, $\mu \not= \nu$,
$\mu \not= 1$, $\nu \not= 1$.
\end{itemize}

\subsection{{\boldmath $\dim \X =4$}}
\label{ss-dim4}

Finally we give the complete list of 4-dimensional unital right coideals of
$\mc{U}$.  For the proof of the completeness of the list the method
explained above is used.

\begin{itemize}
\item
$\X ^4_1=\Lin \{\G 3,\G 2,G,1\}$,
\item
$\X ^4_2=\Lin \{K^2\G 2+\alpha \F 1K^2+\beta K^2\E 1,K^2G,K^2,1\}$,
$|\alpha |+|\beta |\not=0$,
\item
$\X ^4_3=\Lin \{\G 2+\alpha \F 1+\beta \E 1,G,K^{-2},1\}$,
$|\alpha |+|\beta |\not=0$,
\item
$\X ^4_4=\Lin \{f_\mu \G 2,f_\mu G,f_\mu ,1\}$, $\mu \not=1$,
\item
$\X ^4_5=\Lin \{\G 2,G,\F 1K^2+\alpha K^2\E 1+\beta K^2,1\}$,
\item
$\X ^4_6=\Lin \{\G 2,G,K^2\E 1+\alpha K^2,1\}$,
\item
$\X ^4_7=\Lin \{\G 2,G,f_\mu ,1\}$, $\mu \not=1$,
\item
$\X ^4_8=\Lin \{\F 1K^2G+\alpha K^2\E 1G+\beta K^2G+\gamma K^2\E 1
+\delta K^2,G,\\
\F 1K^2+\alpha K^2\E 1+\beta K^2,1\}$,\\
(LI) $\gamma \not=2\alpha $
\item
$\X ^4_9=\Lin \{K^2\E 1G+\alpha K^2G+\beta \F 1K^2+\gamma K^2,G,
K^2\E 1+\alpha K^2,1\}$,\\
(LI) $\beta \not= 0$
\item
$\X ^4_{10}=\Lin \{K^4G+\alpha \F 2K^4+\alpha \beta \F 1K^4\E 1
+\alpha \beta ^2K^4\E 2+\alpha \gamma \F 1K^4+\alpha \beta \gamma K^4\E 1,
\F 1K^2+\beta K^2\E 1+\gamma K^2,K^4,1\}$, $\alpha \not=0$,
\item
$\X ^4_{11}=\Lin \{K^4G+\alpha K^4\E 2+\alpha \beta K^4\E 1,
K^2\E 1+\beta K^2,K^4,1\}$, $\alpha \not=0$,
\item
$\X ^4_{12}=\Lin \{G+\alpha \F 2+\alpha \beta \F 1\E 1
+\alpha \beta ^2\E 2+\alpha \gamma \F 1+\alpha \beta \gamma \E 1,\\
\F 1K^{-2}+\beta K^{-2}\E 1+\gamma K^{-2},K^{-4},1\}$, $\alpha \not=0$,
\item
$\X ^4_{13}=\Lin \{G+\alpha \E 2+\alpha \beta \E 1,
K^{-2}\E 1+\beta K^{-2},K^{-4},1\}$, $\alpha \not=0$,
\item
$\X ^4_{14}=\Lin \{f_\mu G+\alpha \F 1f_\mu +\beta f_\mu \E 1,
f_\mu ,f_{\mu q^{-1}},1\}$,\\
$\mu \not=1$, $\mu \not=q$, $|\alpha |+|\beta |\not=0$,
\item
$\X ^4_{15}=\Lin \{K^2 G+\alpha \F 1K^2 +\beta K^2 \E 1,
G,K^2,1\}$, $|\alpha |+|\beta |\not=0$,
\item
$\X ^4_{16}=\Lin \{K^2G+\alpha \F 1K^2 +\beta K^2\E 1,
\F 1K^4+\gamma K^4\E 1+\delta K^4,K^2,1\}$,
$|\alpha |+|\beta |\not=0$,\\
(LI) $\beta \not= \alpha \gamma $,
\item
$\X ^4_{17}=\Lin \{K^2G+\alpha \F 1K^2 +\beta K^2\E 1,
K^4\E 1+\gamma K^4,K^2,1\}$,\\ $|\alpha |+|\beta |\not=0$,\\
(LI) $\alpha \not= 0$,
\item
$\X ^4_{18}=\Lin \{K^2G+\alpha K^2\E 1,
\F 1K^2+\beta K^2\E 1,K^2,1\}$, $\alpha \not=0$,\\
(LI)
\item
$\X ^4_{19}=\Lin \{K^2G+\alpha \F 1K^2,K^2\E 1,K^2,1\}$,
$\alpha \not=0$,\\
(LI)
\item
$\X ^4_{20}=\Lin \{K^2G+\alpha \F 1K^2+\beta K^2\E 1,f_\mu ,K^2,1\}$,
$|\alpha |+|\beta |\not=0$, $\mu \not=1$, $\mu \not= q$,
\item
$\X ^4_{21}=\Lin \{G+\alpha \F 1 +\beta \E 1,
K^{-2}G,K^{-2},1\}$, $|\alpha |+|\beta |\not=0$,
\item
$\X ^4_{22}=\Lin \{G+\alpha \E 1,\F 1+\beta \E 1,K^{-2},1\}$,
$\alpha \not=0$,\\
(LI)
\item
$\X ^4_{23}=\Lin \{G+\alpha \F 1,\E 1,K^{-2},1\}$, $\alpha \not=0$,\\
(LI)
\item
$\X ^4_{24}=\Lin \{G+\alpha \F 1 +\beta \E 1,
\F 1K^2+\gamma K^2\E 1+\delta K^2,K^{-2},1\}$,\\
$|\alpha |+|\beta |\not=0$,\\
(LI) $\beta \not= \alpha \gamma $
\item
$\X ^4_{25}=\Lin \{G+\alpha \F 1 +\beta \E 1,
K^2\E 1+\gamma K^2,K^{-2},1\}$, $|\alpha |+|\beta |\not=0$,\\
(LI) $\alpha \not= 0$
\item
$\X ^4_{26}=\Lin \{G+\alpha \F 1 +\beta \E 1,
f_\mu ,K^{-2},1\}$, $\mu \not= 1$, $\mu \not= q^{-1}$,
\item
$\X ^4_{27}=\Lin \{f_\mu G,G,f_\mu ,1\}$, $\mu \not= 1$,
\item
$\X ^4_{28}=\Lin \{f_\mu G,\F 1f_{\mu q}+\alpha f_{\mu q}\E 1
+\beta f_{\mu q},f_\mu ,1\}$, $\mu \not= 1$,
\item
$\X ^4_{29}=\Lin \{f_\mu G,f_{\mu q}\E 1+\alpha f_{\mu q},
f_\mu ,1\}$, $\mu \not= 1$,
\item
$\X ^4_{30}=\Lin \{f_\mu G,\F 1K^2+\alpha K^2\E 1+\beta K^2,
f_\mu ,1\}$, $\mu \not= 1$,
\item
$\X ^4_{31}=\Lin \{f_\mu G,K^2\E 1+\alpha K^2,f_\mu ,1\}$,
$\mu \not= 1$,
\item
$\X ^4_{32}=\Lin \{f_\mu G,f_\mu ,f_\nu ,1\}$, $\mu \not= \nu $,
$\mu \not= 1$, $\nu \not= 1$,
\item
$\X ^4_{33}=\Lin \{G,\F 2K^4+\alpha \F 1K^4\E 1+\alpha ^2K^4\E 2
+\beta \F 1K^4+\alpha \beta K^4\E 1+ \gamma K^4,
\F 1K^2+\alpha K^2\E 1+\beta K^2,1\}$, $\alpha \not= 0$,
\item
$\X ^4_{34}=\Lin \{G,K^4\E 2+\alpha K^4\E 1+\beta K^4,
K^2\E 1+\alpha K^2,1\}$,
\item
$\X ^4_{35}=\Lin \{G,\F 1f_\mu +\alpha f_\mu \E 1+\beta f_\mu ,
f_{\mu q^{-1}},1\}$, $\mu \not= q$,
\item
$\X ^4_{36}=\Lin \{G,f_\mu \E 1+\alpha f_\mu ,f_{\mu q^{-1}},1\}$,
$\mu \not= q$,
\item
$\X ^4_{37}=\Lin \{G,\F 1K^2+\alpha K^2, K^2\E 1+\beta K^2,1\}$,\\
(LI)
\item
$\X ^4_{38}=\Lin \{G,\F 1K^2+\alpha K^2\E 1+\beta K^2, f_\mu ,1\}$,
$\mu \not= 1$,
\item
$\X ^4_{39}=\Lin \{G,K^2\E 1+\alpha K^2, f_\mu ,1\}$, $\mu \not= 1$,
\item
$\X ^4_{40}=\Lin \{G,f_\mu ,f_\nu ,1\}$, $\mu \not= \nu$,
$\mu \not= 1$, $\nu \not= 1$,
\item
$\X ^4_{41}=\Lin \{\F 3K^6+\alpha \F 2K^6\E 1+\alpha ^2\F 1K^6\E 2
+\alpha ^3K^6\E 3+\\
\beta \F 2K^6+\alpha \beta \F 1K^6\E 1+\alpha ^2\beta K^6\E 2
+\gamma \F 1K^6 +\alpha \gamma K^6\E 1+\delta K^6,\\
\F 2K^4+\alpha \F 1K^4\E 1+\alpha ^2K^4\E 2+\beta \F 1K^4+\alpha \beta K^4\E 1
+\gamma K^4,\\
\F 1K^2+\alpha K^2\E 1+\beta K^2,1\}$,
\item
$\X ^4_{42}=\Lin \{K^6\E 3+\alpha K^6\E 2+\beta K^6\E 1+\gamma K^6,
K^4\E 2+\alpha K^4\E 1+\beta K^4,K^2\E 1+\alpha K^2,1\}$,
\item
$\X ^4_{43}=\Lin \{\F 2f_\mu +\alpha \F 1f_\mu \E 1+\alpha ^2 f_\mu
\E 2+\beta \F 1f_\mu +\alpha \beta f_\mu \E 1+\gamma f_\mu ,
\F 1f_{\mu q^{-1}}+\alpha f_{\mu q^{-1}}\E 1+\beta f_{\mu q^{-1}},
f_{\mu q^{-2}},1\}$, $\mu \not= q^2$,
\item
$\X ^4_{44}=\Lin \{f_\mu \E 2+\alpha f_\mu \E 1+\beta f_\mu ,
f_{\mu q^{-1}}\E 1+\alpha f_{\mu q^{-1}},f_{\mu q^{-2}},1\}$, $\mu \not= q^2$,
\item
$\X ^4_{45}=\Lin \{\F 2K^4+\alpha _1\F 1K^4\E 1
+\alpha _2K^4\E 2+(\beta _1+\alpha _1\beta _2)\F 1K^4
+(\alpha _1\beta _1+\alpha _2\beta _2)K^4\E 1+\gamma K^4,
\F 1K^2+\beta _1K^2,K^2\E 1+\beta _2K^2,1\}$, $\alpha _2\not= \alpha _1^2$,\\
(LI) $(q^2-q^{-2})\gamma \not= q\alpha _1+(q^2-1)(\beta _1^2+2\alpha _1\beta _1
\beta _2+\alpha _2\beta _2^2)$
\item
$\X ^4_{46}=\Lin \{\F 2K^4+\alpha \F 1K^4\E 1
+\alpha ^2K^4\E 2+\beta \F 1K^4+\gamma K^4\E 1+\delta K^4,
\F 1K^2+\alpha K^2\E 1,K^2,1\}$, $\gamma \not= \alpha \beta $,\\
(LI)
\item
$\X ^4_{47}=\Lin \{\F 1K^4\E 1+\alpha K^4\E 2+\beta \F 1K^4
+(\gamma +\alpha \beta )K^4\E 1+\delta K^4,
\F 1K^2+\gamma K^2,K^2\E 1+\beta K^2,1\}$,\\
(LI) $(q^2-q^{-2})\delta \not= q+(q^2-1)\beta (2\gamma +\alpha \beta )$
\item
$\X ^4_{48}=\Lin \{K^4\E 2+\alpha \F 1K^4 +\beta K^4\E 1+\gamma K^4,
K^2\E 1,K^2,1\}$, $\alpha \not= 0$,\\
(LI)
\item
$\X ^4_{49}=\Lin \{\F 2K^4+\alpha \F 1K^4\E 1+\alpha ^2K^4\E 2
+\beta \F 1K^4+\alpha \beta K^4\E 1+\gamma K^4,
\F 1K^2+(\beta-\alpha \delta)K^2, K^2\E 1+\delta K^2,1\}$,\\
(LI) $\alpha \not= (q-q^{-3})\gamma -(q-q^{-1})\beta ^2$
\item
$\X ^4_{50}=\Lin \{\F 2K^4+\alpha \F 1K^4\E 1+\alpha ^2K^4\E 2
+\beta \F 1K^4+\alpha \beta K^4\E 1+\gamma K^4,
\F 1K^2+\alpha K^2\E 1+\beta K^2, f_\mu ,1\}$, $\mu \not= 1$,
\item
$\X ^4_{51}=\Lin \{K^4\E 2+\alpha K^4\E 1+\beta K^4,
\F 1K^2+\gamma K^2,K^2\E 1+\alpha K^2,1\}$,\\
(LI) $(1+q^{-2})\beta \not= \alpha ^2$
\item
$\X ^4_{52}=\Lin \{K^4\E 2+\alpha K^4\E 1+\beta K^4,
K^2\E 1+\alpha K^2,f_\mu,1\}$, $\mu \not= 1$,
\item
$\X ^4_{53}=\Lin \{\F 1f_\mu +\alpha f_\mu \E 1+\beta f_\mu ,
\F 1K^2+\gamma K^2\E 1+\delta K^2,f_{\mu q^{-1}},1\}$, $\mu \not= q$,\\
(LI) $\gamma \not= \alpha $, $\mu ^2\not= q^2$
\item
$\X ^4_{54}=\Lin \{\F 1f_\mu +\alpha f_\mu \E 1+\beta f_\mu ,
K^2\E 1+\gamma K^2,f_{\mu q^{-1}},1\}$, $\mu \not= q$,\\
(LI) $\mu ^2\not= q^2$
\item
$\X ^4_{55}=\Lin \{\F 1f_\mu +\alpha f_\mu ,f_\mu \E 1+\beta f_\mu ,
f_{\mu q^{-1}},1\}$, $\mu \not= q$,\\
(LI) $\mu ^2\not= q^2$
\item
$\X ^4_{56}=\Lin \{\F 1f_\mu +\alpha f_\mu \E 1+\beta f_\mu,
f_{\mu q^{-1}},f_\nu ,1\}$, $\mu \not= q$, $\mu \not= q\nu $, $\nu \not= 1$,
\item
$\X ^4_{57}=\Lin \{\F 1K^2+\alpha K^2\E 1+\beta K^2,
f_\mu \E 1+\gamma f_\mu ,f_{\mu q^{-1}},1\}$, $\mu \not= q$,\\
(LI) $\mu ^2\not= q^2$
\item
$\X ^4_{58}=\Lin \{\F 1K^2+\alpha K^2,K^2\E 1+\beta K^2,f_\mu ,1\}$,
$\mu \not= 1$,\\
(LI) $\mu ^2\not= 1$
\item
$\X ^4_{59}=\Lin \{\F 1K^2+\alpha K^2\E 1+\beta K^2,f_\mu ,f_\nu ,
1\}$, $\mu \not= \nu $, $\mu \not= 1$, $\nu \not= 1$,
\item
$\X ^4_{60}=\Lin \{K^2\E 1+\alpha K^2,f_\mu ,f_\nu ,1\}$,
$\mu \not= \nu $, $\mu \not= 1$, $\nu \not= 1$,
\item
$\X ^4_{61}=\Lin \{f_\mu \E 1+\alpha f_\mu ,K^2\E 1+\beta K^2,
f_{\mu q^{-1}},1\}$, $\mu \not= q$,
\item
$\X ^4_{62}=\Lin \{f_\mu \E 1+\alpha f_\mu ,f_{\mu q^{-1}},f_\nu ,
1\}$, $\mu \not= q$, $\mu \not= q\nu $, $\nu \not= 1$,
\item
$\X ^4_{63}=\Lin \{f_{\mu _1},f_{\mu _2},f_{\mu _3},1\}$,
$\mu _i\not= \mu _j$, $\mu _i\not= 1$ for any $i\not= j$.
\end{itemize}

\section{Further restrictions on the calculus}
\label{sec-restr}

\subsection{Generators of the FODC}

Let $\Ga $ be a 3-dimensional left-covariant differential calculus on
$\SLq 2$. We shall look for calculi $\Ga $ satisfying the following additional
condition:
\begin{itemize}
\item[(LI)]
The left-invariant one-forms $\omega (u^1_2)$, $\omega (u^2_1)$ and
$\omega (u^1_1-u^2_2)$ generate $\Ga $ as a left $\A $-module.
\end{itemize}
It is simple to check that this condition is fulfilled if and only if
the subspace $\bar{\X }$ of the linear functionals on the matrix elements of
the fundamental corepresentation $u$ is four dimensional.
Note that the classical differential calculus on $\mathrm{SL}(2)$ obviously
has this property.
%
All coideals from Subsection~\ref{ss-dim4} satisfying condition (LI) are marked
with the label (LI).
The necessary and sufficient conditions for the parameter values
are indicated after the label.

\subsection{Universal (higher order) differential calculi}
\label{ss-univcal}

In the previous section we have seen that there is a very large number of
left-covariant first order differential calculi on the quantum group
$\SLq 2$. Now the corresponding left-covariant universal (higher order)
differential calculi will be considered, too. We require that
\begin{itemize}
\item
the dimension of the space of left-invariant differential 2-forms is
at least 3.
\end{itemize}
Recall that any differential calculus is a quotient of the universal
differential calculus. Hence
if there is a calculus such that the dimension of the space of its
left-invariant differential 2-forms is equal to 3 then the universal calculus
satisfies the above condition.

In order to study this condition we use Lemma \ref{l-dimGaw}
and the list in Subsection \ref{ss-dim4}.
As a sample let us consider $\X =\X ^4_8$, $\gamma \not= 2\alpha $.
We then have
\begin{align*}
\bar{\X }&=\{\x _1,\x _2,\x _3,\x _4\}\\
&=\{FKG+\alpha KEG+\beta K^2G+\gamma KE+\delta K^2,G,FK+\alpha KE+\beta K^2,1\}.
\end{align*}
For the elements $\mathrm{m}(\x _i\ot \x _j)$ we obtain
\allowdisplaybreaks
\begin{align}
\label{eq-x1x1}
\x _1\x _1=&q^{-1}F^2K^2G^2+(1+q^{-2})\alpha FK^2EG^2
+q^{-1}\alpha ^2 K^2E^2G^2\\
\notag &+\text{terms of lower degree}\\
\label{eq-x1x2}
\x _1\x _2=& FKG^2+\alpha KEG^2+\beta K^2G^2+\gamma KEG+\delta K^2G\\
\label{eq-x1x3}
\x _1\x _3=& q^{-1}F^2K^2G+(1+q^{-2})\alpha FK^2EG+q^{-1}\alpha ^2K^2E^2G\\
\notag &+\text{terms of lower degree}\\
\label{eq-x1x4}
\x _1\x _4=& FKG+\alpha KEG+\beta K^2G+\gamma KE+\delta K^2\\
\label{eq-x2x1}
\x _2\x _1-&\x _1\x _2+2\x _1\x _4=
4\alpha KEG+2\beta K^2G+4\gamma KE+2\delta K^2\\
\label{eq-x2x2}
\x _2\x _2=& G^2\\
\label{eq-x2x3}
\x _2\x _3-&\x _1\x _4+2\x _3\x _4=
(-\gamma +4\alpha )KE+(-\delta +2\beta )K^2\\
\label{eq-x2x4}
\x _2\x _4=&G\\
\label{eq-x3x1}
\x _3\x _1-&\x _1\x _3-2\x _3^2=
((1-q^{-2})\gamma -4\alpha )FK^2E-4q^{-1}\alpha ^2K^2E^2\\
\notag &+\text{terms of lower degree}\\
\label{eq-x3x2}
\x _3\x _2-&\x _1\x _4= -\gamma KE-\delta K^2\\
\label{eq-x3x3}
\x _3\x _3=&
q^{-1}F^2K^2+(1+q^{-2})\alpha FK^2E+q^{-1}\alpha ^2K^2E^2
+(1+q^{-2})\beta FK^3 \notag \\
&+(1+q^{-2})\alpha \beta K^3E
+(\beta ^2+\alpha /(q-q^{-1}))K^4-\alpha /(q-q^{-1})\\
\label{eq-x3x4}
\x _3\x _4=&FK+\alpha KE+\beta K^2\\
\label{eq-x4x1}
\x _4\x _1-&\x _1\x _4=0\\
\label{eq-x4x2}
\x _4\x _2-&\x _2\x _4=0\\
\label{eq-x4x3}
\x _4\x _3-&\x _3\x _4=0\\
\label{eq-x4x4}
\x _4\x _4=&1
\end{align}
{}From Lemma \ref{l-dimGaw} we conclude that there are at least
$3+\dim \X =6$
linearly independent relations in $\bar{\X }^2$. Hence there are at most
10 linearly independent elements in $\bar{\X }^2$.
Because of the PBW-theorem the 9 elements in (\ref{eq-x1x1})--(\ref{eq-x1x4}),
(\ref{eq-x2x2}), (\ref{eq-x2x4}), (\ref{eq-x3x3}), (\ref{eq-x3x4}) and
(\ref{eq-x4x4}) are linearly independent.
Suppose that $\alpha $ is nonzero. Then (\ref{eq-x2x1}) and the difference
(\ref{eq-x2x3})$-$(\ref{eq-x3x2}) are further linearly independent elements
which is a contradiction. If $\alpha =0$ then $\gamma \not= 0$ because of
(LI). Therefore, (\ref{eq-x2x3}) and (\ref{eq-x3x1}) are two additional
linearly independent elements in $\bar{\X }^2$ and we obtain again a
contradiction.

The same procedure can be applied to all right coideals of the list in
Subsection \ref{ss-dim4}. In order to carry out these computations we used
the computer algebra program FELIX \cite{inp-ApelKlaus}.
The result is the following (the number of the general coideal
and the corresponding parameter values are given in parenthesis):
\begin{itemize}
\item[1.]
(47; $\alpha =\beta =\gamma =0$, $\delta=q^5/(q^2-1)^2$)\\
$\bar{\X }=\{FK^2E+q^5/(q^2-1)^2 K^4,FK,KE,1\}$
\item[2.]
(53; $\alpha =0$, $\beta \not=0$, $\gamma =-(q-q^{-1})^2\beta ^2$,
$\delta=(1+q)\beta $, $\mu =q^{1/2}$)\\
$\bar{\X }=\{F+\beta K,FK-(q-q^{-1})^2\beta ^2KE+(1+q)\beta K^2,K^{-1},1\}$
\item[3.]
(53; $\alpha =0$, $\beta \not=0$, $\gamma =-(q-q^{-1})^2\beta ^2$,
$\delta=(1+q)\beta $, $\mu =-q^{1/2}$)\\
$\bar{\X }=\{F\coun _-+\beta \coun _-K,
FK-(q-q^{-1})^2\beta ^2KE+(1+q)\beta K^2,\coun _-K^{-1},1\}$
\item[4.]
(53; $\alpha =0$, $\beta \not=0$, $\gamma =(q-q^{-1})^2\beta ^2$,
$\delta=(1-q)\beta $, $\mu =\ii q^{1/2}$)\\
$\bar{\X }=\{Ff_\ii +\beta f_\ii K,
FK +(q-q^{-1})^2\beta ^2KE+(1-q)\beta K^2,f_\ii K^{-1},1\}$
\item[5.]
(53; $\alpha =0$, $\beta \not=0$, $\gamma =(q-q^{-1})^2\beta ^2$,
$\delta=(1-q)\beta $, $\mu =-\ii q^{1/2}$)\\
$\bar{\X }=\{Ff_{-\ii }+\beta f_{-\ii }K,
FK +(q-q^{-1})^2\beta ^2KE+(1-q)\beta K^2,f_{-\ii }K^{-1},1\}$
\item[6.]
(53; $\alpha \not= 0$, $\gamma =-\alpha $, $\beta =\delta =0$, $\mu =\ii q$)
\qquad $\dim \uGa ^{\wedge 2}=4$\\
$\bar{\X }=\{Ff_\ii K+\alpha f_\ii KE, FK -\alpha KE, f_\ii , 1\}$
\item[7.]
(53; $\alpha \not= 0$, $\gamma =-\alpha $, $\beta =\delta =0$, $\mu =-\ii q$)
\qquad $\dim \uGa ^{\wedge 2}=4$\\
$\bar{\X }=\{Ff_{-\ii }K+\alpha f_{-\ii }KE, FK -\alpha KE, f_{-\ii }, 1\}$
\item[8.]
(54; $\alpha =\beta =\gamma =0$, $\mu =q^3$)\\
$\bar{\X }=\{FK^5,KE,K^4,1\}$
\item[9.]
(54; $\alpha =\beta =\gamma =0$, $\mu =q^{-1}$)\\
$\bar{\X }=\{FK^{-3},KE,K^{-4},1\}$
\item[10.]
(55; $\alpha =\beta =0$, $\mu =1$)\\
$\bar{\X }=\{FK^{-1},K^{-1}E,K^{-2},1\}$
\item[11.]
(55; $\alpha =\beta =0$, $\mu =-1$)\\
$\bar{\X }=\{F\coun _-K^{-1},\coun _-K^{-1}E,\coun _-K^{-2},1\}$
\item[12.]
(55; $\alpha =\beta =0$, $\mu =q^{-1}$)\\
$\bar{\X }=\{FK^{-3},K^{-3}E,K^{-4},1\}$
\item[13.]
(57; $\alpha =\beta =\gamma =0$, $\mu =q^3$)\\
$\bar{\X }=\{FK,K^5E,K^4,1\}$
\item[14.]
(57; $\alpha =\beta =\gamma =0$, $\mu =q^{-1}$)\\
$\bar{\X }=\{FK,K^{-3}E,K^{-4},1\}$
\item[15.]
(57; $\gamma \not= 0$, $\alpha =-1/((q-q^{-1})^2\gamma ^2)$,
$\beta =1/((q-1)(q^{-2}-1)\gamma )$, $\mu =q^{1/2}$)\\
$\bar{\X }=\{FK+\alpha KE+\beta K^2, E+\gamma K, K^{-1}, 1\}$
\item[16.]
(57; $\gamma \not= 0$, $\alpha =-1/((q-q^{-1})^2\gamma ^2)$,
$\beta =1/((q-1)(q^{-2}-1)\gamma )$, $\mu =-q^{1/2}$)\\
$\bar{\X }=\{FK+\alpha KE+\beta K^2, \coun _-E+\gamma \coun _-K,
\coun _-K^{-1}, 1\}$
\item[17.]
(58; $\alpha =\beta =0$, $\mu =q^2$)\\
$\bar{\X }=\{FK,KE,K^4,1\}$
\item[18.]
(58; $\alpha =\beta =0$, $\mu =q$)\\
$\bar{\X }=\{FK,KE,K^2,1\}$
\item[19.]
(58; $\alpha \not= 0$, $\beta =-q^3/((q^2-1)^2\alpha )$,
$\mu =q^{-1}$)\\
$\bar{\X }=\{FK+\alpha K^2,KE-q^3/((q^2-1)^2\alpha )K^2,K^{-2},1\}$
\item[20.]
(58; $\alpha =\beta =0$, $\mu =-q$)\\
$\bar{\X }=\{FK,KE,\coun _-K^2,1\}$
\end{itemize}

\begin{bem}
The solutions 6 and 7 have the property
$\dim \uGa ^{\wedge 2}=4$. In all
other cases we have $\dim \uGa ^{\wedge 2}=3$.
\end{bem}

\subsection{Hopf algebra automorphisms}

Let $\A $ be a Hopf algebra and $\varphi $ an automorphism of $\A $.
Let $\Ga $ be an $\A $-bimodule. Then the actions
$.:\A \times \Ga \to \Ga $, $a.\rho :=\varphi (a)\rho $ and 
$.:\Ga \times \A \to \Ga $, $\rho .a:=\rho \varphi (a)$ determine
another $\A $-bimodule structure on $\Ga $, since $\varphi $ is an
algebra homomorphism. Let $\Ga $ be a left $\A $-module. Then the mapping
$\lcoa ':\Ga \to \A \ot \Ga $, $\lcoa '(\rho )=\varphi (\rho _{(-1)}) \ot
\rho _{(0)}$ determines a second left $\A $-comodule structure on $\Ga $,
since $\varphi $ is a comodule homomorphism.

If $(\Ga ,\dif )$ is a left-covariant FODC over $\A $, set
$\Ga _\varphi :=\Ga $ with module action $.$ and comodule mapping $\lcoa '$.
Further define $\dif _\varphi :\A \to \Ga $,
$\dif _\varphi a:=\dif \varphi (a)$. Then
$(\Ga _\varphi ,\dif _\varphi )$ is a left-covariant FODC over $\A $. Indeed,
$\dif _\varphi $ satisfies the Leibniz rule and
\[ \Lin \{a.\dif _\varphi b\,|\,a,b\in \A \}=
\Lin \{\varphi (a)\dif \varphi (b)\,|\,a,b\in \A \}=
\Lin \{a\dif b\,|\,a,b\in \A \}\]
since $\varphi $ is invertible.

\begin{defin}
Let $(\Ga ,\dif )$ be a left-covariant FODC over the Hopf algebra $\A $.
We call $(\Ga ,\dif )$ \textit{Hopf-invariant} if
$(\Ga _\varphi ,\dif _\varphi )$ is isomorphic to $(\Ga ,\dif )$ for all
Hopf algebra automorphism $\varphi $ of $\A $.
\end{defin}

\begin{satz}
Let $(\Ga ,\dif )$ be a left-covariant FODC over $\A $. The following
statements are equivalent:\\
(i) $(\Ga ,\dif )$ is Hopf-invariant.\\
(ii) $\varphi (\mc{R}_\Ga )=\mc{R}_\Ga $ for any Hopf algebra automorphism
$\varphi $ of $\A $.\\
(iii) $\varphi (\X _\Ga )=\X _\Ga $ for any Hopf algebra automorphism
$\varphi $ of $\Anull $ which is implemented by $\A $.
\end{satz}

\begin{bew}
It suffices to show that for any Hopf algebra automorphism $\varphi $ of $\A $
we have $\mc{R}_{\Ga _\varphi }=\varphi ^{-1}(\mc{R}_\Ga )$
and $\X _{\Ga _\varphi }=\varphi (\X _\Ga )$.
For this we compute
\begin{align*}
\omega _\varphi (a):=S(a_{(1)}).\dif _\varphi a_{(2)}
=\varphi (S(a_{(1)}))\dif \varphi (a_{(2)})
=S(\varphi (a)_{(1)})\dif \varphi (a)_{(2)}
=\omega (\varphi (a))
\end{align*}
for any $a\in \A $. Hence (with $\A ^+=\A \cap \ker \coun $)
\begin{align*}
\mc{R}_{\Ga _\varphi }&=\{a\in \A ^+\,|\,\omega _\varphi (a)=0\}\\
&=\{a\in \A ^+\,|\,\omega (\varphi (a))=0\}
=\{\varphi ^{-1}(b)\in \A ^+\,|\,\omega (b)=0\}
=\varphi ^{-1}(\mc{R}_\Ga )
\end{align*}
and (with $\Anull {}^+:=\Anull \cap \ker \coun $)
\begin{align*}
\X _{\Ga _\varphi }&=\{X\in \Anull {}^+\,|\,X(\mc{R}_{\Ga _\varphi })=0\}
=\{X\in \Anull {}^+\,|\,X(\varphi ^{-1}(\mc{R}_\Ga ))=0\}\\
&=\{Y\circ \varphi\in \Anull {}^+\,|\,Y(\mc{R}_\Ga )=0\}
=\varphi (\X _\Ga ).
\end{align*}
\end{bew}

It is easy to prove that any Hopf algebra automorphism $\varphi $ of
$\OSLq 2$ is given by
$\varphi (u^1_1)=u^1_1$,
$\varphi (u^2_2)=u^2_2$,
$\varphi (u^1_2)=\alpha u^1_2$,
$\varphi (u^2_1)=\alpha ^{-1}u^2_1$, where
$\alpha \in \compx $ (see also \cite[Section 4.1.2]{b-KS}). For this
$\varphi =:\varphi _\alpha $ we obtain from (\ref{eq-pairing}) the formulas
\begin{align}
\varphi _\alpha (E)=\alpha ^{-1}E,\quad \varphi _\alpha (F)=\alpha F,\quad
\varphi _\alpha (f_\mu )=f_\mu ,\qquad \mu \in \compx .
\end{align}

\begin{satz}
A left-covariant FODC $\Ga $ over $\OSLq 2$ is
Hopf-invariant if and only if its quantum tangent space $\X _\Ga $ is
generated by elements which are homogeneous with respect to the
$\mathbb{Z}$-grading of $\mc{U}$.
\end{satz}

\begin{bew}
Clearly $\varphi _\alpha (X)=\alpha ^{-n}X$ for any $X\in \X $,
$\mathrm{deg}X=n$. Therefore the condition of the Proposition is sufficient.

Conversely, suppose that $X=\sum _{i=1}^k\lambda _iX_i\in \X $, $k>1$, where
$\lambda _i\in \compx $, $X_i\in \X $, $\mathrm{deg}X_i=n_i$
and $n_i\not= n_j$ for any $i\not= j$. If $(\Ga ,\dif )$ is Hopf-invariant
then we have $\X \ni \varphi _q(X)-q^{-n_k}X=
\sum _{i=1}^{k-1}(q^{-n_i}-q^{-n_k})X_i$. Since $q$ is not a root of unity,
we obtain that $\sum _{i=1}^{k-1}\mu _iX_i\in \X $, $\mu _i\in \compx $.
By induction on $k$ one easily proves that $X_i\in \X $ for any $i$.
\end{bew}

\begin{folg}\label{f-Hopfinvcal}
The items 1, 8-14, 17, 18 and 20 are precisely the quantum tangent spaces of
the list in Subsection \ref{ss-univcal} which are Hopf-invariant.

\end{folg}

\section{More structures on left-covariant differential calculi on $\SLq 2$}
\label{sec-struc}

The Hopf algebra $\OSLq 2$ admits 3 non-equivalent real forms. Namely,
\begin{itemize}
\item
$q\in \mathbb{R}$;
$(u^1_1)^*=u^2_2$, $(u^1_2)^*=-qu^2_1$, $(u^2_1)^*=-q^{-1}u^1_2$,
$(u^2_2)^*=u^1_1$;
\item
$q\in \mathbb{R}$;
$(u^1_1)^*=u^2_2$, $(u^1_2)^*=qu^2_1$, $(u^2_1)^*=q^{-1}u^1_2$,
$(u^2_2)^*=u^1_1$;
\item
$|q|=1$;
$(u^i_j)^*=u^i_j$ for any $i,j=1,2$.
\end{itemize}
The Hopf $*$-algebras corresponding to them are
$\mc{O}(\mathrm{SU}_q(2))$, $\mc{O}(\mathrm{SU}_q(1,1))$ and
$\mc{O}(\mathrm{SL}_q(2,\mathbb{R}))$.
Let us introduce the dual involution on $\mc{U}$ in the way
$f^*(a):=\overline{f(S(a)^*)}$. Then the corresponding $*$-structures on
$\mc{U}$ are given by
\begin{itemize}
\item
$E^*=F$, $F^*=E$, $G^*=G$, $f_\mu ^*=f_{\bar{\mu }}$,
\item
$E^*=-F$, $F^*=-E$, $G^*=G$, $f_\mu ^*=f_{\bar{\mu }}$,
\item
$E^*=-qE$, $F^*=-q^{-1}F$, $G^*=-G$, $f_\mu ^*=f_{\bar{\mu }^{-1}}$,
\end{itemize}
respectively. A FODC $(\Ga ,\dif )$ over a Hopf $*$-algebra $\A $ is called a
$*$-\textit{calculus} if there exists an involution $*:\Ga \to \Ga $ such that
$(a\,(\dif b)\,c)^*=c^*(\dif b^*)a^*$ for any $a,b,c\in \A $.

\begin{satz}\cite{b-KS}
Let $(\Ga ,\dif )$ be a finite-dimensional left-covariant FODC over a Hopf
$*$-algebra $\A $. Then $(\Ga ,\dif )$ is a $*$-calculus if and only if
its quantum tangent space $\X _\Ga $ is $*$-invariant.
\end{satz}

Let $\Ga $ be a left-covariant bimodule over a Hopf algebra $\A $. An
invertible linear mapping $\sigma :\Ga \otA \Ga \to \Ga \otA \Ga $ is called
a \textit{braiding of} $\Ga $ if $\sigma $
is a homomorphism of $\A $-bimodules, commutes with the left coaction on $\Ga $
and satisfies the braid relation
\[
(\sigma \ot \id )(\id \ot \sigma )(\sigma \ot \id )=
(\id \ot \sigma )(\sigma \ot \id )(\id \ot \sigma )
\]
on $\Ga ^{\ot 3}$.
If $(\Ga ,\dif )$ is a left-covariant differential calculus over $\A $ then
we require that $(\id -\sigma )(\mc{S}\cap \Ga ^{\ot 2})=0$.
Such a braiding neither needs to exist nor it is unique for a given
left-covariant differential calculus over $\A $.

Let $(\Ga ,\dif )$ be a left-covariant FODC over $\A $. Fix a basis
$\{\x _i\,|\,i=1,\ldots ,n\}$ of its quantum tangent space $\X _\Ga $
and let $\{\w _i\,|\,i=1,\ldots ,n\}$ be the dual basis of $\linv{\Ga }$.
If $\sigma $ is a braiding of $\Ga $ with $\sigma (\w _i\ot \w _j)=
\sigma ^{kl}_{ij}\w _k\ot \w _l$ then we can define a bilinear mapping
$[\cdot ,\cdot ]:\X _\Ga \times \X _\Ga \to \X _\Ga $ by
\begin{align}
[\x _i,\x _j]:=\x _i\x _j-\sigma ^{ij}_{kl}\x _k\x _l.
\end{align}
This mapping can be viewed as a generalization of the Lie bracket of
the left-invariant vector fields if it also satisfies a generalized
Jacobi identity. For the calculi and braidings described below
the mapping $\beta :=[\cdot ,\cdot ]$ fulfills the equation
\begin{align}
(\beta \circ (\beta \ot \id ) -\beta \circ (\id \ot \beta ))A_3^\mathrm{t}=0
\end{align}
where $A_3^\mathrm{t}(\x _i\ot \x _j \ot \x _k)=
(\id -\sigma _{12})(\id -\sigma _{23}+\sigma _{23}\sigma _{12})^{ijk}_{rst}
\x _r\ot \x _s\ot \x _t$.
Unfortunately, such an equation does not hold for bicovariant
differential calculi in general.

\section{The calculi in detail}
\label{sec-detail}

The main result of the considerations of the preceding sections is the
following.

\begin{thm}\label{t-main}
Suppose that $q$ is a transcendental complex number. Let
$(\Ga ,\dif )$ be a left-covariant FODC over $\OSLq 2$ having the
following properties:
\begin{itemize}
\item
$\Ga $ is a 3-dimensional left-covariant bimodule
\item
the left-invariant one-forms $\omega (u^1_2),\omega (u^2_1)$ and
$\omega (u^1_1-u^2_2)$ form a basis of the left module $\Ga $
\item
the universal differential calculus $\uGa ^\wedge $ associated to $\Ga $
satisfies the inequality $\dim \uGa ^{\wedge 2}\geq 3$
\item
$\Ga $ is Hopf-invariant.
\end{itemize}
Then $\Ga $ is isomorphic to one of the 11 left-covariant FODC
1,8-14,17,18,20 of the list in Subsection \ref{ss-univcal}.
\end{thm}

\begin{bems}
1. The calculus 17 is isomorphic to the 3D-calculus of Woronowicz
\cite{a-Woro3}. Two other 3-dimensional calculi appear in \cite{a-SchSch3}.
They are isomorphic to the calculus 13 (for $r=2$) and 8 (for $r=3$),
respectively. The calculi 10 and 11 are subcalculi of the bicovariant
$4\mathrm{D}_+$- and $4\mathrm{D}_-$-calculus, respectively (see also
\cite[Section 14.2.4]{b-KS}). One can easily check all these isomorphisms
by comparing the corresponding right ideals $\mc{R}_\Ga $.
\\
2. In \cite[Section 7]{a-Schmue1} the notion of \textit{elementary}
left-covariant
FODC was introduced. It is an easy computation to show that the six
calculi 10--12,17,18 and 20 in Theorem \ref{t-main} are elementary, while
the others are not.
\end{bems}

After the last section we list some important facts about these 11 calculi.
First we fix the basis
$\wH :=\omega ((u^1_1-u^2_2)/2),\wX :=\omega (u^1_2),\wY :=\omega (u^2_1)$
in $\linv{\Ga }$ and determine the dual basis $\{H,X,Y\}$
(i.\,e.\ $\omega (a)=H(a)\wH +X(a)\wX +Y(a)\wY $ for any $a\in \OSLq 2$).
We examine whether or not the calculus is a $*$-calculus with respect to the
three given
involutions on $\OSLq 2$. We compute the pairing between the quantum tangent
space and linear (fundamental representation) and quadratic elements of
$\OSLq 2$. Using this the generators of the right ideal $\mc{R}_\Ga $
are computed.

The quadratic-linear relations between generators of $\X _\Ga $ are given.
Because of their simple form one can easily show that the algebra without unit
generated by the elements $H,X$, and $Y$ and these relations has the PBW-basis
$\{H^{n_1}X^{n_2}Y^{n_3}\,|\,n_1,n_2,n_3\in \mathbb{N}_0, n_1+n_2+n_3>0\}$.
This proves also that $\dim \uGa ^{\wedge 3}=1$ and $\dim \uGa ^{\wedge k}=0$
for any $k>3$ and for each of the 11 calculi.

The matrix $(f^i_j)$ describes the commutation rules between one-forms
and functions (see equation (\ref{eq-module})), where $i,j\in \{H,X,Y\}$.
Finally, the generators of the right ideal and equation (\ref{eq-symmform})
determine the vector space of left-invariant symmetric 2-forms.
If there exists a braiding $\sigma $ of the left-covariant bimodule $\Ga $ such
that $(\id -\sigma )(\mc{S}\cap \linv{\Ga }^{\ot 2})=0$ then we give one of
them by its eigenvalues and eigenspaces. It was determined by means of
the computer algebra program FELIX \cite{inp-ApelKlaus}.

\section{Cohomology}
\label{sec-deRham}

Let us fix one of the differential calculi of Theorem \ref{t-main}.
The differential mapping $\dif :\Ga ^\wedge \to \Ga ^\wedge $ satisfies
the equation $\dif ^2=0$. Hence it defines a complex
\begin{gather}\label{eq-dcomplex}
\begin{CD}
\{0\}@>\dif >>\A @>\dif >>\Ga @>\dif >>\Ga ^{\wedge 2}
@>\dif >>\Ga ^{\wedge 3} @>\dif >>\Ga ^{\wedge 4}=\{0\}.
\end{CD}
\end{gather}
Now the corresponding cohomology spaces will be determined.
In \cite{a-Woro3} Woro\-nowicz introduced a method to do this. Let us recall
his strategy.

The first observation is that the algebra $\A =\OSLq 2$ has a vector space
basis
\begin{gather}
\{v^\lambda _{ij}\,|\,\lambda \in \frac{1}{2}\mathbb{N}_0;
-\lambda \leq i,j\leq \lambda ;\lambda -i,\lambda -j\in \mathbb{N}_0\}
\end{gather}
such that $\copr (v^\lambda _{ij})=v^\lambda _{ik}\ot v^\lambda _{kj}$
for any $\lambda ,i,j$. Let us fix such a basis. Then for
$\rho _j\in \linv{\Ga }^\wedge $ we obtain
\begin{gather}\label{eq-difform}
\dif (v^\lambda _{ij}\rho _j )=
v^\lambda _{ik}\x _r(v^\lambda _{kj})\w _r\wedge \rho _j
+v^\lambda _{ij}\dif \rho _j.
\end{gather}
The elements $\sum _{r,j}\x _r(v^\lambda _{kj})\w _r\wedge \rho _j
+\dif \rho _k$ are left-invariant. Hence the differential mapping $\dif $
preserves the direct sum decomposition
\begin{gather}
\Ga ^\wedge =\bigoplus _{\lambda ,i}C(v^\lambda _i)\linv{\Ga }^\wedge ,
\end{gather}
where $C(v^\lambda _i)=\Lin \{v^\mu _{kl}\,|\,\mu =\lambda ,k=i\}$.
Moreover, formula (\ref{eq-difform}) for the differential on
$C(v^\lambda _i)\linv{\Ga }^\wedge $ does not depend on $i$.
Therefore, for any $\lambda \in \frac{1}{2}\mathbb{N}_0$ we obtain $2\lambda $
isomorphic differential complexes.

Let $V^\lambda =\Lin \{e^\lambda _\mu \,|\,-\lambda \leq \mu \leq \lambda ;
\lambda -\mu \in \mathbb{N}_0\}$ be the representation of $\Anull $
defined by $f.e^\lambda _\mu =e^\lambda _\nu f(v^\lambda _{\nu \mu })$
for any $f\in \Anull $. It is isomorphic to the representation given by
\begin{align}\label{eq-anullrep}
E.e^\lambda _\nu &=[\lambda +\nu ]e^\lambda _{\nu -1},&
F.e^\lambda _\nu &=[\lambda -\nu ]e^\lambda _{\nu +1},\notag \\
f_\mu .e^\lambda _\nu &=\mu ^{-2\nu }e^\lambda _\nu ,&
G.e^\lambda _\nu &=-2\nu e^\lambda _\nu .
\end{align}
Because of the left-covariance of the differential calculus, $V^\lambda $
induces a representation of the quantum tangent space
$\mc{X}_\Ga $ of the differential calculus. We obtain a new complex
\begin{gather}\label{eq-deRhamseq}
\begin{CD}
\{0\}@>\dif >>V^\lambda @>\dif >>V^\lambda \ot \linv{\Ga }@>\dif >>
V^\lambda \ot \linv{\Ga }^{\wedge 2}\\
@. @. @>\dif >>V^\lambda \ot \linv{\Ga }^{\wedge 3}@>\dif >>\{0\},
\end{CD}
\end{gather}
where
\begin{gather}\label{eq-deRham}
\dif (e^\lambda _\mu \ot \rho _\mu )=
e^\lambda _\nu \x _r(v^\lambda _{\nu \mu })\ot \w _r\wedge \rho _\mu
+e^\lambda _\mu \ot \dif \rho _\mu .
\end{gather}

\begin{lemma}\label{l-deRham}
The cohomology spaces of the complex (\ref{eq-deRham}) are isomorphic to
$H^0_\lambda =H^3_\lambda =\comp ^p$,
$H^1_\lambda =H^2_\lambda =\{0\}$,
where $p=1$ for $\lambda =0$ and $p=0$ otherwise.
\end{lemma}

\begin{thm}\label{t-cohom}
Let $\Ga $ be one of the differential calculi of Theorem \ref{t-main}.
Then the cohomology spaces of the differential complex (\ref{eq-dcomplex})
are isomorphic to
\begin{align}
H^0&=H^3=\comp ^p,&
H^1&=H^2=\{0\},
\end{align}
where $p=1$ for $\lambda =0$ and $p=0$ otherwise.
\end{thm}

\begin{bew}
The assertion follows from Lemma \ref{l-deRham} and the preceding
considerations.
\end{bew}

\begin{bew}[ of the Lemma]
For $\lambda =0$ we have $\dif e^0_0=0$ and $\dif (e^0_0\ot \w )\not= 0$
for any $\w \in \linv{\Ga }$.
Since $\dim \linv{\Ga }^{\wedge 2}=\dim \linv{\Ga }(=3)<\infty $, the mapping
$\dif :\linv{\Ga }\to \linv{\Ga }^{\wedge 2}$ is also surjective. Therefore,
for any $\xi \in \linv{\Ga }^{\wedge 2}$ we obtain $\dif \xi (=\dif ^2\w )=0$.
Hence the assertion of the lemma for $\lambda =0$ is valid.

Let now $\lambda \not= 0$. We set
\begin{gather}
\partial (e^\lambda _\mu \ot \bigwedge _{j=1}^k\w _{i_j}):=\mu
+\sum _{j=1}^k\partial (\w _{i_j}),\text{ where }
\partial (\wH )=0,
\partial (\wY )=-\partial (\wX )=1.
\end{gather}
Then $\partial $ defines a direct sum decomposition of
$V^\lambda \ot \linv{\Ga }^\wedge $ and $\dif $ preserves this decomposition:
$\partial (\dif \rho )=\partial (\rho )$ for any homogeneous element
$\rho $ of $V^\lambda \ot \linv{\Ga }^\wedge $.
Hence it suffices to check that the sequence (\ref{eq-deRhamseq}) restricted
to homogeneous elements of the same degree $\mu $ is exact for any
$\lambda \in \frac{1}{2}\mathbb{N}$ and any $\mu $ with
$\lambda -\mu \in \mathbb{Z}$.

For $\rho =e^\lambda _\mu $ the equation
$\dif \rho (=\x _i.e^\lambda _\mu \ot \w _i)=0$ implies that
$X.e^\lambda _\mu =Y.e^\lambda _\mu =0$. Hence $\mu =-\lambda $ and
$\mu =\lambda $. This is a contradiction to $\lambda \not= 0$, so we have
$H^0_\lambda =\{0\}$.

Since $\dif \rho =0$ for any $\rho \in V^\lambda \ot \linv{\Ga }^{\wedge 3}$
we have to show that $\rho =\dif \xi $ for some
$\xi \in V^\lambda \ot \linv{\Ga }^{\wedge 2}$.
If $\partial (\rho )=\mu $ then $-\lambda \leq \mu \leq \lambda $ and
$\rho $ is a constant multiple of
$e^\lambda _\mu \ot \wH \wedge \wX \wedge \wY $. One can take
$\xi _1:=c^{-1}e^\lambda _{\mu +1}\ot \wH \wedge \wX $ or
$\xi _2:=c^{-1}e^\lambda _{\mu -1}\ot \wH \wedge \wY $ with some
$c\in \compx $. Indeed,
because of $\wH \wedge \wH \wedge \wX =\wX \wedge \wH \wedge \wX =
\dif (\wH \wedge \wX )=0$ the element
$\dif (e^\lambda _{\mu +1}\ot \wH \wedge \wX )=Y.e^\lambda _{\mu +1}\ot
\wY \wedge \wH \wedge \wX $ is a nonzero multiple of $\rho $ for
$\mu \not=\lambda $.
Similarly, $\dif (c\xi _2)$ is a nonzero multiple
of $\rho $ for $\mu \not= -\lambda $. Therefore, $H^3_\lambda =\{0\}$.

Now we prove that $H^1_\lambda =\{0\}$. It is easy to see that if
$\partial \w =\mu $ for a nonzero $\w \in V^\lambda \ot \linv{\Ga }$ then
$|\mu |\leq \lambda +1$. Hence we have to show that
$\dif \w \not= 0$ if $|\mu |=\lambda +1$ and that the vector space
$\{\w \in V^\lambda \ot \linv{\Ga }\,|\,\partial \w =\mu ,\dif \w =0\}$
is one-dimensional for $|\mu |\leq \lambda $.
We define $V^{\lambda \mu }_1:=\{\w \in V^\lambda \ot \linv{\Ga }\,|\,
\partial \w =\mu \}$.
Since $e^\lambda _\lambda \ot \wY $ generates the vector space
$V^{\lambda ,\lambda +1}_1$ we have to show that
$H.e^\lambda _\lambda \ot \wH \wedge \wY +e^\lambda _\lambda \ot \dif \wY
\not= 0$. Similarly,
$\dif (e^\lambda _{-\lambda }\ot \wX)=
H.e^\lambda _{-\lambda }\ot \wH \wedge \wX
+e^\lambda _{-\lambda }\ot \dif \wX $
must be nonzero. Using the explicit formulas for the quantum tangent space,
for the differentials $\dif \wX $ and $\dif \wY $ and (\ref{eq-anullrep})
this is easily done.

Secondly, the elements $e^\lambda _\mu \ot \wH $,
$e^\lambda _{\mu +1}\ot \wX $ and $e^\lambda _{\mu -1}\ot \wY $ generate
the vector space $V^{\lambda \mu }_1$. We have $\dim V^{\lambda \mu }_1=2$
for $|\mu |=\lambda $ and $\dim V^{\lambda \mu }_1=3$ for $|\mu |<\lambda $.
Moreover,
\begin{align}
\notag
\dif (e^\lambda _\mu \ot \wH )=&
X.e^\lambda _\mu \ot \wX \wedge \wH
+Y.e^\lambda _\mu \ot \wY \wedge \wH\\
\label{eq-dH}
&+(H.e^\lambda _\mu \ot \wH \wedge \wH
+e^\lambda _\mu \ot \dif \wH )\\
\label{eq-dX}
\dif (e^\lambda _{\mu +1}\ot \wX )=&
Y.e^\lambda _{\mu +1}\ot \wY \wedge \wX
+(H.e^\lambda _{\mu +1}\ot \wH \wedge \wX
+e^\lambda _{\mu +1}\ot \dif \wX )
\\
\label{eq-dY}
\dif (e^\lambda _{\mu -1}\ot \wY )=&
X.e^\lambda _{\mu -1}\ot \wX \wedge \wY
+(H.e^\lambda _{\mu -1}\ot \wH \wedge \wY
+e^\lambda _{\mu -1}\ot \dif \wY )
\end{align}
Observe that $\xi = \cHX (\xi )\ot \wH \wedge \wX
+\cHY (\xi )\ot \wH \wedge \wY +\cXY (\xi )\ot \wX \wedge \wY $
for any $\xi \in V^\lambda \ot \linv{\Ga }^{\wedge 2}$, where
$\cHX (\xi ),\cHY (\xi ),\cXY (\xi )\in V^\lambda $.
Now if $\mu =\lambda $ then by (\ref{eq-dY})
$\cXY (\dif (e^\lambda _{\lambda -1}\ot \wY ))=X.e^\lambda _{\lambda -1}
\not= 0$
and therefore the range of $\dif $ is at least one-dimensional.
Similarly, by (\ref{eq-dX})
$\cXY (\dif (e^\lambda _{\mu +1}\ot \wX ))\not= 0$ for $\mu =-\lambda $.
Hence $\dim \ker \dif \!\!\upharpoonright \!\!V^{\lambda \mu }_1\leq 1$ for
$|\mu |=\lambda $.
If $|\mu |<\lambda $ then
$\cHY (\dif (e^\lambda _\mu \ot \wH ))\not= 0$,
$\cHY (\dif (e^\lambda _{\mu +1}\ot \wX ))= 0$ and
$\cXY (\dif (e^\lambda _{\mu +1}\ot \wX ))\not= 0$.
Therefore, $\dim \dif (V^{\lambda \mu }_1)\geq 2$ and so
$\dim \ker \dif \!\!\upharpoonright \!\!V^{\lambda \mu }_1\leq 1$.
Together,
$\dim \ker \dif \!\!\upharpoonright \!\!V^{\lambda \mu }_1\leq 1$
for $|\mu |\leq \lambda $ and
$\ker \dif \!\!\upharpoonright \!\!V^{\lambda \mu }_1=\{0\}$ for
$|\mu |=\lambda +1$.
Because of $\dif e^\lambda _\mu \not= 0$,
$\partial (\dif e^\lambda _\mu )=\mu $ and
$\dif (\dif e^\lambda _\mu )=0$ we also have
$\dim \ker \dif \!\!\upharpoonright \!\!V^{\lambda \mu }_1\geq 1 $ for
$|\mu |\leq \lambda $.
This means that $H^1_\lambda =\{0\}$.

The last assertion, $H^2_\lambda =\{0\}$, follows from dimension
computations.
\end{bew}

\newpage

\noindent
\textbf{1.\ Quantum tangent space {\boldmath $\X _\Ga $}:}
\begin{align*}
H&:=\frac{2(q^{-1}-q)}{q^4+1}\left( FK^2E+\frac{q^5(K^4-1)}{(q^2-1)^2}\right),
& X&:=q^{-1/2}FK, & Y&:=q^{-1/2}KE
\end{align*}

\noindent
\textbf{Real forms:}
$\mc{O}(\mathrm{SU}_q(2))$,
$\mc{O}(\mathrm{SU}_q(1,1))$,
$\mc{O}(\mathrm{SL}_q(2,\mathbb{R}))$

\noindent
\textbf{Fund. repr.:}
$\displaystyle H=\frac{2}{q^2+q^{-2}}
\begin{pmatrix}q^{-2} & 0\\0 & -q^2\end{pmatrix}$,
$X=\begin{pmatrix}0 & 1\\0 & 0\end{pmatrix}$,
$Y=\begin{pmatrix}0 & 0\\1 & 0\end{pmatrix}$

\noindent
\textbf{Relations:}
\begin{align*}
qHX-q^{-1}XH&=\frac{2(q+q^{-1})}{q^2+q^{-2}}X\\
q^{-1}HY-qYH&=\frac{-2(q+q^{-1})}{q^2+q^{-2}}Y\\
q^{-2}XY-q^2YX&=\frac{q^2+q^{-2}}{2}H
\end{align*}

\noindent
\textbf{Module structure, differentials:}
\begin{align*}
(f^i_j)=
\begin{pmatrix}
K^4 & 0 & 0\\
\frac{2q^{-3/2}(q^{-1}-q)}{q^2+q^{-2}}K^3E & K^2 & 0\\
\frac{2q^{-3/2}(q^{-1}-q)}{q^2+q^{-2}}FK^3 & 0 & K^2
\end{pmatrix}
\quad &
\begin{array}{rl}
\dif \wH =& -\frac{q^4+1}{2}\wX \wedge \wY \\
\dif \wX =& \frac{-2(q^{-2}+1)}{q^2+q^{-2}}\wH \wedge \wX \\
\dif \wY =& \frac{2(q^2+1)}{q^2+q^{-2}}\wH \wedge \wY
\end{array}
\end{align*}

\noindent
\textbf{Pairing:}
\[
\begin{array}{r|ccccccccc}
 & (u^1_1)^2 & u^1_1u^1_2 & u^1_1u^2_1 & (u^1_2)^2 & u^1_2u^2_1 & u^1_2u^2_2
& (u^2_1)^2 & u^2_1u^2_2 & (u^2_2)^2\\
\hline
H & \frac{2q^{-2}+2q^{-4}}{q^2+q^{-2}} & 0 & 0 & 0 &
\frac{2(q^{-1}-q)}{q^2+q^{-2}} & 0 & 0 & 0 & \frac{-2q^4-2q^2}{q^2+q^{-2}}\\
X & 0 & 1 & 0 & 0 & 0 & q & 0 & 0 & 0\\
Y & 0 & 0 & 1 & 0 & 0 & 0 & 0 & q & 0
\end{array}
\]

\noindent
\textbf{Right ideal} {\boldmath $\mc{R}_\Ga $}\textbf{:}
$u^1_1+q^{-4}u^2_2-(1+q^{-4})$, $(u^1_2)^2$, $(u^2_1)^2$,
$u^1_2u^2_1+(q^3-q)u^1_1$, $(u^1_1-1)u^1_2$, $(u^1_1-1)u^2_1$

\noindent
\textbf{Left-invariant symmetric 2-forms:}
\begin{align*}
&\wH \ot \wH
&
&\wX \ot \wX
&
&q^{-1}\wH \ot \wX +q\wX \ot \wH
\\
&q^2\wX \ot \wY +q^{-2}\wY \ot \wX
&
&\wY \ot \wY
&
&q\wH \ot \wY +q^{-1}\wY \ot \wH
\end{align*}

\noindent
\textbf{Braiding:} $(1-\sigma )(q^2+\sigma )=0$\\
$\ker (q^2+\sigma):$
$q^{-2}\wH \ot \wX -q^2\wX \ot \wH $,
$q^2\wH \ot \wY -q^{-2}\wY \ot \wH $,\\
\phantom{$\ker (q^2+\sigma):$}
$q^3\wX \ot \wY -q^{-3}\wY \ot \wX
 -\frac{4(q-q^{-1})}{(q^2+q^{-2})^2}\wH \ot \wH $.

\newpage

\noindent
\textbf{2.\ Quantum tangent space {\boldmath $\X _\Ga $}:}
\begin{align*}
H&:=\frac{2}{q^{-2}-q^2}( K^4-1), & X&:=q^{-5/2}FK^5, & Y&:=q^{-1/2}KE
\end{align*}

\noindent
\textbf{Real forms:}
$\mc{O}(\mathrm{SL}_q(2,\mathbb{R}))$

\noindent
\textbf{Fund. repr.:}
$\displaystyle H=\frac{2}{q+q^{-1}}
\begin{pmatrix}q^{-1} & 0\\0 & -q\end{pmatrix}$,
$X=\begin{pmatrix}0 & 1\\0 & 0\end{pmatrix}$,
$Y=\begin{pmatrix}0 & 0\\1 & 0\end{pmatrix}$

\noindent
\textbf{Relations:}
\begin{align*}
q^2HX-q^{-2}XH&=2X\\
q^{-2}HY-q^2YH&=-2Y\\
q^{-3}XY-q^3YX+\frac{(q^2+1)^2(q^2-1)}{4q^3}H^2&=\frac{q+q^{-1}}{2}H
\end{align*}

\noindent
\textbf{Module structure, differentials:}
\begin{align*}
(f^i_j)=
\begin{pmatrix}
K^4 & \frac{q^2-q^{-2}}{2}X & 0\\
0 & K^6 & 0\\
0 & 0 & K^2
\end{pmatrix}
\quad &
\begin{array}{rl}
\dif \wH =& \frac{-q^4-q^2}{2}\wX \wedge \wY \\
\dif \wX =& -2q^{-2}\wH \wedge \wX \\
\dif \wY =& 2q^2\wH \wedge \wY
\end{array}
\end{align*}

\noindent
\textbf{Pairing:}
\[
\begin{array}{r|ccccccccc}
 & (u^1_1)^2 & u^1_1u^1_2 & u^1_1u^2_1 & (u^1_2)^2 & u^1_2u^2_1 & u^1_2u^2_2
& (u^2_1)^2 & u^2_1u^2_2 & (u^2_2)^2\\
\hline
H & 2q^{-2} & 0 & 0 & 0 & 0 & 0 & 0 & 0 & -2q^2\\
X & 0 & q^{-2} & 0 & 0 & 0 & q^3 & 0 & 0 & 0\\
Y & 0 & 0 & 1 & 0 & 0 & 0 & 0 & q & 0
\end{array}
\]

\noindent
\textbf{Right ideal} {\boldmath $\mc{R}_\Ga $}\textbf{:}
$u^1_1+q^{-2}u^2_2-(1+q^{-2})$, $(u^1_2)^2$, $(u^2_1)^2$,
$u^1_2u^2_1$, $(u^1_1-q^{-2})u^1_2$, $(u^1_1-1)u^2_1$

\noindent
\textbf{Left-invariant symmetric 2-forms:} $p=(q^2+1)^2(q^2-1)/4$
\begin{align*}
& \wH \ot \wH -p\wX \ot \wY
&
&\wX \ot \wX
&
&q^{-2}\wH \ot \wX +q^2\wX \ot \wH
\\
&q^3\wX \ot \wY +q^{-3}\wY \ot \wX
&
&\wY \ot \wY
&
&q^2\wH \ot \wY +q^{-2}\wY \ot \wH
\end{align*}

\noindent
\textbf{Braiding:} ---

\newpage

\noindent
\textbf{3.\ Quantum tangent space {\boldmath $\X _\Ga $}:}
\begin{align*}
H&:=\frac{2}{q^2-q^{-2}}( K^{-4}-1), & X&:=q^{3/2}FK^{-3}, & Y&:=q^{-1/2}KE
\end{align*}

\noindent
\textbf{Real forms:}
$\mc{O}(\mathrm{SL}_q(2,\mathbb{R}))$

\noindent
\textbf{Fund. repr.:}
$\displaystyle H=\frac{2}{q+q^{-1}}
\begin{pmatrix}q & 0\\0 & -q^{-1}\end{pmatrix}$,
$X=\begin{pmatrix}0 & 1\\0 & 0\end{pmatrix}$,
$Y=\begin{pmatrix}0 & 0\\1 & 0\end{pmatrix}$

\noindent
\textbf{Relations:}
\begin{align*}
q^{-2}HX-q^2XH&=2X\\
q^2HY-q^{-2}YH&=-2Y\\
qXY-q^{-1}YX&=\frac{q+q^{-1}}{2}H
\end{align*}

\noindent
\textbf{Module structure, differentials:}
\begin{align*}
(f^i_j)=
\begin{pmatrix}
K^{-4} & \frac{q^2-q^{-2}}{2}X & 0\\
0 & K^{-2} & 0\\
0 & 0 & K^2
\end{pmatrix}
\quad &
\begin{array}{rl}
\dif \wH =& \frac{-1-q^{-2}}{2}\wX \wedge \wY \\
\dif \wX =& -2q^2\wH \wedge \wX \\
\dif \wY =& 2q^{-2}\wH \wedge \wY
\end{array}
\end{align*}

\noindent
\textbf{Pairing:}
\[
\begin{array}{r|ccccccccc}
 & (u^1_1)^2 & u^1_1u^1_2 & u^1_1u^2_1 & (u^1_2)^2 & u^1_2u^2_1 & u^1_2u^2_2
& (u^2_1)^2 & u^2_1u^2_2 & (u^2_2)^2\\
\hline
H & 2q^2 & 0 & 0 & 0 & 0 & 0 & 0 & 0 & -2q^{-2}\\
X & 0 & q^2 & 0 & 0 & 0 & q^{-1} & 0 & 0 & 0\\
Y & 0 & 0 & 1 & 0 & 0 & 0 & 0 & q & 0
\end{array}
\]

\noindent
\textbf{Right ideal} {\boldmath $\mc{R}_\Ga $}\textbf{:}
$u^1_1+q^2u^2_2-(1+q^2)$, $(u^1_2)^2$, $(u^2_1)^2$,
$u^1_2u^2_1$, $(u^1_1-q^2)u^1_2$, $(u^1_1-1)u^2_1$

\noindent
\textbf{Left-invariant symmetric 2-forms:}
\begin{align*}
&\wH \ot \wH
&
&\wX \ot \wX
&
&q^2\wH \ot \wX +q^{-2}\wX \ot \wH
\\
&q^{-1}\wX \ot \wY +q\wY \ot \wX
&
&\wY \ot \wY
&
&q^{-2}\wH \ot \wY +q^2\wY \ot \wH
\end{align*}

\noindent
\textbf{Braiding:} $(1-\sigma )(q^2+\sigma )=0$\\
$\ker (q^2+\sigma):$
$q\wH \ot \wX -q^{-1}\wX \ot \wH $,
$q^{-1}\wH \ot \wY -q\wY \ot \wH $,\\
\phantom{$\ker (q^2+\sigma):$}
$\wX \ot \wY -\wY \ot \wX $.

\newpage

\noindent
\textbf{4.\ Quantum tangent space {\boldmath $\X _\Ga $}:}
\begin{align*}
H&:=\frac{2}{q-q^{-1}}( K^{-2}-1), & X&:=q^{1/2}FK^{-1}, & Y&:=q^{1/2}K^{-1}E
\end{align*}

\noindent
\textbf{Real forms:}
$\mc{O}(\mathrm{SU}_q(2))$,
$\mc{O}(\mathrm{SU}_q(1,1))$,
$\mc{O}(\mathrm{SL}_q(2,\mathbb{R}))$

\noindent
\textbf{Fund. repr.:}
$\displaystyle H=\frac{2}{q+1}
\begin{pmatrix}q & 0\\0 & -1\end{pmatrix}$,
$X=\begin{pmatrix}0 & 1\\0 & 0\end{pmatrix}$,
$Y=\begin{pmatrix}0 & 0\\1 & 0\end{pmatrix}$

\noindent
\textbf{Relations:}
\begin{align*}
q^{-1}HX-qXH&=2X\\
qHY-q^{-1}YH&=-2Y\\
qXY-q^{-1}YX-\frac{q-q^{-1}}{4}H^2&=H
\end{align*}

\noindent
\textbf{Module structure, differentials:}
\begin{align*}
(f^i_j)=
\begin{pmatrix}
K^{-2} & \frac{q-q^{-1}}{2}X & \frac{q-q^{-1}}{2}Y\\
0 & 1 & 0\\
0 & 0 & 1
\end{pmatrix}
\quad &
\begin{array}{rl}
\dif \wH =& -q^{-1}\wX \wedge \wY \\
\dif \wX =& -2q\wH \wedge \wX \\
\dif \wY =& 2q^{-1}\wH \wedge \wY
\end{array}
\end{align*}

\noindent
\textbf{Pairing:}
\[
\begin{array}{r|ccccccccc}
 & (u^1_1)^2 & u^1_1u^1_2 & u^1_1u^2_1 & (u^1_2)^2 & u^1_2u^2_1 & u^1_2u^2_2
& (u^2_1)^2 & u^2_1u^2_2 & (u^2_2)^2\\
\hline
H & 2q & 0 & 0 & 0 & 0 & 0 & 0 & 0 & -2q^{-1}\\
X & 0 & q & 0 & 0 & 0 & 1 & 0 & 0 & 0\\
Y & 0 & 0 & q & 0 & 0 & 0 & 0 & 1 & 0
\end{array}
\]

\noindent
\textbf{Right ideal} {\boldmath $\mc{R}_\Ga $}\textbf{:}
$u^1_1+qu^2_2-(1+q)$, $(u^1_2)^2$, $(u^2_1)^2$,
$u^1_2u^2_1$, $(u^1_1-q)u^1_2$, $(u^1_1-q)u^2_1$

\noindent
\textbf{Left-invariant symmetric 2-forms:}
\begin{align*}
&\wH \ot \wH +\frac{1-q^{-2}}{4}\wX \ot \wY
&
&\wX \ot \wX
&
&q\wH \ot \wX +q^{-1}\wX \ot \wH
\\
&q^{-1}\wX \ot \wY +q\wY \ot \wX
&
&\wY \ot \wY
&
&q^{-1}\wH \ot \wY +q\wY \ot \wH
\end{align*}

\noindent
\textbf{Braiding:} $(1-\sigma )(q^2+\sigma )=0$\\
$\ker (q^2+\sigma):$
$\wH \ot \wX -\wX \ot \wH $,
$\wH \ot \wY -\wY \ot \wH $,
$\wX \ot \wY -\wY \ot \wX $.

\newpage

\noindent
\textbf{5.\ Quantum tangent space {\boldmath $\X _\Ga $}:}
\begin{align*}
H&:=\frac{2}{q^{-1}-q}(\coun _-K^{-2}-1), & X&:=-q^{1/2}F\coun _-K^{-1}, &
Y&:=-q^{1/2}\coun _-K^{-1}E
\end{align*}

\noindent
\textbf{Real forms:}
$\mc{O}(\mathrm{SU}_q(2))$,
$\mc{O}(\mathrm{SU}_q(1,1))$,
$\mc{O}(\mathrm{SL}_q(2,\mathbb{R}))$

\noindent
\textbf{Fund. repr.:}
$\displaystyle H=\frac{2}{q-1}
\begin{pmatrix}q & 0\\0 & 1\end{pmatrix}$,
$X=\begin{pmatrix}0 & 1\\0 & 0\end{pmatrix}$,
$Y=\begin{pmatrix}0 & 0\\1 & 0\end{pmatrix}$

\noindent
\textbf{Relations:}
\begin{align*}
q^{-1}HX-qXH&=-2X\\
qHY-q^{-1}YH&=2Y\\
qXY-q^{-1}YX-\frac{q-q^{-1}}{4}H^2&=-H
\end{align*}

\noindent
\textbf{Module structure, differentials:}
\begin{align*}
(f^i_j)=
\begin{pmatrix}
\coun _-K^{-2} & \frac{q^{-1}-q}{2}X & \frac{q^{-1}-q}{2}Y\\
0 & \coun _- & 0\\
0 & 0 & \coun _-
\end{pmatrix}
\quad &
\begin{array}{rl}
\dif \wH =& q^{-1}\wX \wedge \wY \\
\dif \wX =& 2q\wH \wedge \wX \\
\dif \wY =& -2q^{-1}\wH \wedge \wY
\end{array}
\end{align*}

\noindent
\textbf{Pairing:}
\[
\begin{array}{r|ccccccccc}
 & (u^1_1)^2 & u^1_1u^1_2 & u^1_1u^2_1 & (u^1_2)^2 & u^1_2u^2_1 & u^1_2u^2_2
& (u^2_1)^2 & u^2_1u^2_2 & (u^2_2)^2\\
\hline
H & -2q & 0 & 0 & 0 & 0 & 0 & 0 & 0 & 2q^{-1}\\
X & 0 & -q & 0 & 0 & 0 & -1 & 0 & 0 & 0\\
Y & 0 & 0 & -q & 0 & 0 & 0 & 0 & -1 & 0
\end{array}
\]

\noindent
\textbf{Right ideal} {\boldmath $\mc{R}_\Ga $}\textbf{:}
$u^1_1-qu^2_2-(1-q)$, $(u^1_2)^2$, $(u^2_1)^2$,
$u^1_2u^2_1$, $(u^1_1+q)u^1_2$, $(u^1_1+q)u^2_1$

\noindent
\textbf{Left-invariant symmetric 2-forms:}
\begin{align*}
&\wH \ot \wH +\frac{1-q^{-2}}{4}\wX \ot \wY
&
&\wX \ot \wX
&
&q\wH \ot \wX +q^{-1}\wX \ot \wH
\\
&q^{-1}\wX \ot \wY +q\wY \ot \wX
&
&\wY \ot \wY
&
&q^{-1}\wH \ot \wY +q\wY \ot \wH
\end{align*}

\noindent
\textbf{Braiding:} $(1-\sigma )(q^2+\sigma )=0$\\
$\ker (q^2+\sigma):$
$\wH \ot \wX -\wX \ot \wH $,
$\wH \ot \wY -\wY \ot \wH $,
$\wX \ot \wY -\wY \ot \wX $.

\newpage

\noindent
\textbf{6.\ Quantum tangent space {\boldmath $\X _\Ga $}:}
\begin{align*}
H&:=\frac{2}{q^2-q^{-2}}( K^{-4}-1), & X&:=q^{3/2}FK^{-3}, & Y&:=q^{3/2}K^{-3}E
\end{align*}

\noindent
\textbf{Real forms:}
$\mc{O}(\mathrm{SU}_q(2))$,
$\mc{O}(\mathrm{SU}_q(1,1))$,
$\mc{O}(\mathrm{SL}_q(2,\mathbb{R}))$

\noindent
\textbf{Fund. repr.:}
$\displaystyle H=\frac{2}{q+q^{-1}}
\begin{pmatrix}q & 0\\0 & -q^{-1}\end{pmatrix}$,
$X=\begin{pmatrix}0 & 1\\0 & 0\end{pmatrix}$,
$Y=\begin{pmatrix}0 & 0\\1 & 0\end{pmatrix}$

\noindent
\textbf{Relations:}
\begin{align*}
q^{-2}HX-q^2XH&=2X\\
q^2HY-q^{-2}YH&=-2Y\\
q^3XY-q^{-3}YX-\frac{(q^2+1)^2(q^2-1)}{4q^3}H^2&=\frac{q+q^{-1}}{2}H
\end{align*}

\noindent
\textbf{Module structure, differentials:}
\begin{align*}
(f^i_j)=
\begin{pmatrix}
K^{-4} & \frac{q^2-q^{-2}}{2}X & \frac{q^2-q^{-2}}{2}Y\\
0 & K^{-2} & 0\\
0 & 0 & K^{-2}
\end{pmatrix}
\quad &
\begin{array}{rl}
\dif \wH =& \frac{-q^{-2}-q^{-4}}{2}\wX \wedge \wY \\
\dif \wX =& -2q^2\wH \wedge \wX \\
\dif \wY =& 2q^{-2}\wH \wedge \wY
\end{array}
\end{align*}

\noindent
\textbf{Pairing:}
\[
\begin{array}{r|ccccccccc}
 & (u^1_1)^2 & u^1_1u^1_2 & u^1_1u^2_1 & (u^1_2)^2 & u^1_2u^2_1 & u^1_2u^2_2
& (u^2_1)^2 & u^2_1u^2_2 & (u^2_2)^2\\
\hline
H & 2q^2 & 0 & 0 & 0 & 0 & 0 & 0 & 0 & -2q^{-2}\\
X & 0 & q^2 & 0 & 0 & 0 & q^{-1} & 0 & 0 & 0\\
Y & 0 & 0 & q^2 & 0 & 0 & 0 & 0 & q^{-1} & 0
\end{array}
\]

\noindent
\textbf{Right ideal} {\boldmath $\mc{R}_\Ga $}\textbf{:}
$u^1_1+q^2u^2_2-(1+q^2)$, $(u^1_2)^2$, $(u^2_1)^2$,
$u^1_2u^2_1$, $(u^1_1-q^2)u^1_2$, $(u^1_1-q^2)u^2_1$

\noindent
\textbf{Left-invariant symmetric 2-forms:} $p=(1+q^{-2})^2(1-q^{-2})/4$
\begin{align*}
&\wH \ot \wH +p\wX \ot \wY
&
&\wX \ot \wX
&
&q^2\wH \ot \wX +q^{-2}\wX \ot \wH
\\
&q^{-3}\wX \ot \wY +q^3\wY \ot \wX
&
&\wY \ot \wY
&
&q^{-2}\wH \ot \wY +q^2\wY \ot \wH
\end{align*}

\noindent
\textbf{Braiding:} $(1-\sigma )(q^2+\sigma )=0$\\
$\ker (q^2+\sigma):$
$q\wH \ot \wX -q^{-1}\wX \ot \wH $,
$q^{-1}\wH \ot \wY -q\wY \ot \wH $,\\
\phantom{$\ker (q^2+\sigma):$}
$q^{-2}\wX \ot \wY -q^2\wY \ot \wX $.

\newpage

\noindent
\textbf{7.\ Quantum tangent space {\boldmath $\X _\Ga $}:}
\begin{align*}
H&:=\frac{2}{q^{-2}-q^2}( K^4-1), & X&:=q^{-1/2}FK, & Y&:=q^{-5/2}K^5E
\end{align*}

\noindent
\textbf{Real forms:}
$\mc{O}(\mathrm{SL}_q(2,\mathbb{R}))$

\noindent
\textbf{Fund. repr.:}
$\displaystyle H=\frac{2}{q+q^{-1}}
\begin{pmatrix}q^{-1} & 0\\0 & -q\end{pmatrix}$,
$X=\begin{pmatrix}0 & 1\\0 & 0\end{pmatrix}$,
$Y=\begin{pmatrix}0 & 0\\1 & 0\end{pmatrix}$

\noindent
\textbf{Relations:}
\begin{align*}
q^2HX-q^{-2}XH&=2X\\
q^{-2}HY-q^2YH&=-2Y\\
q^{-3}XY-q^3YX+\frac{(q^2+1)^2(q^2-1)}{4q^3}H^2&=\frac{q+q^{-1}}{2}H
\end{align*}

\noindent
\textbf{Module structure, differentials:}
\begin{align*}
(f^i_j)=
\begin{pmatrix}
K^4 & 0 & \frac{q^{-2}-q^2}{2}Y\\
0 & K^2 & 0\\
0 & 0 & K^6
\end{pmatrix}
\quad &
\begin{array}{rl}
\dif \wH =& \frac{-q^2-q^4}{2}\wX \wedge \wY \\
\dif \wX =& -2q^{-2}\wH \wedge \wX \\
\dif \wY =& 2q^2\wH \wedge \wY
\end{array}
\end{align*}

\noindent
\textbf{Pairing:}
\[
\begin{array}{r|ccccccccc}
 & (u^1_1)^2 & u^1_1u^1_2 & u^1_1u^2_1 & (u^1_2)^2 & u^1_2u^2_1 & u^1_2u^2_2
& (u^2_1)^2 & u^2_1u^2_2 & (u^2_2)^2\\
\hline
H & 2q^{-2} & 0 & 0 & 0 & 0 & 0 & 0 & 0 & -2q^2\\
X & 0 & 1 & 0 & 0 & 0 & q & 0 & 0 & 0\\
Y & 0 & 0 & q^{-2} & 0 & 0 & 0 & 0 & q^3 & 0
\end{array}
\]

\noindent
\textbf{Right ideal} {\boldmath $\mc{R}_\Ga $}\textbf{:}
$u^1_1+q^{-2}u^2_2-(1+q^{-2})$, $(u^1_2)^2$, $(u^2_1)^2$,
$u^1_2u^2_1$, $(u^1_1-1)u^1_2$, $(u^1_1-q^{-2})u^2_1$

\noindent
\textbf{Left-invariant symmetric 2-forms:} $p=(q^2+1)^2(q^2-1)/4$
\begin{align*}
&\wH \ot \wH -p\wX \ot \wY
&
&\wX \ot \wX
&
&q^{-2}\wH \ot \wX +q^2\wX \ot \wH
\\
&q^3\wX \ot \wY +q^{-3}\wY \ot \wX
&
&\wY \ot \wY
&
&q^2\wH \ot \wY +q^{-2}\wY \ot \wH
\end{align*}

\noindent
\textbf{Braiding:} ---

\newpage

\noindent
\textbf{8.\ Quantum tangent space {\boldmath $\X _\Ga $}:}
\begin{align*}
H&:=\frac{2}{q^2-q^{-2}}( K^{-4}-1), & X&:=q^{-1/2}FK, & Y&:=q^{3/2}K^{-3}E
\end{align*}

\noindent
\textbf{Real forms:}
$\mc{O}(\mathrm{SL}_q(2,\mathbb{R}))$

\noindent
\textbf{Fund. repr.:}
$\displaystyle H=\frac{2}{q+q^{-1}}
\begin{pmatrix}q & 0\\0 & -q^{-1}\end{pmatrix}$,
$X=\begin{pmatrix}0 & 1\\0 & 0\end{pmatrix}$,
$Y=\begin{pmatrix}0 & 0\\1 & 0\end{pmatrix}$

\noindent
\textbf{Relations:}
\begin{align*}
q^{-2}HX-q^2XH&=2X\\
q^2HY-q^{-2}YH&=-2Y\\
qXY-q^{-1}YX&=\frac{q+q^{-1}}{2}H
\end{align*}

\noindent
\textbf{Module structure, differentials:}
\begin{align*}
(f^i_j)=
\begin{pmatrix}
K^{-4} & 0 & \frac{q^2-q^{-2}}{2}Y\\
0 & K^2 & 0\\
0 & 0 & K^{-2}
\end{pmatrix}
\quad &
\begin{array}{rl}
\dif \wH =& \frac{-1-q^{-2}}{2}\wX \wedge \wY \\
\dif \wX =& -2q^2\wH \wedge \wX \\
\dif \wY =& 2q^{-2}\wH \wedge \wY
\end{array}
\end{align*}

\noindent
\textbf{Pairing:}
\[
\begin{array}{r|ccccccccc}
 & (u^1_1)^2 & u^1_1u^1_2 & u^1_1u^2_1 & (u^1_2)^2 & u^1_2u^2_1 & u^1_2u^2_2
& (u^2_1)^2 & u^2_1u^2_2 & (u^2_2)^2\\
\hline
H & 2q^2 & 0 & 0 & 0 & 0 & 0 & 0 & 0 & -2q^{-2}\\
X & 0 & 1 & 0 & 0 & 0 & q & 0 & 0 & 0\\
Y & 0 & 0 & q^2 & 0 & 0 & 0 & 0 & q^{-1} & 0
\end{array}
\]

\noindent
\textbf{Right ideal} {\boldmath $\mc{R}_\Ga $}\textbf{:}
$u^1_1+q^2u^2_2-(1+q^2)$, $(u^1_2)^2$, $(u^2_1)^2$,
$u^1_2u^2_1$, $(u^1_1-1)u^1_2$, $(u^1_1-q^2)u^2_1$

\noindent
\textbf{Left-invariant symmetric 2-forms:}
\begin{align*}
&\wH \ot \wH
&
&\wX \ot \wX
&
&q^2\wH \ot \wX +q^{-2}\wX \ot \wH
\\
&q^{-1}\wX \ot \wY +q\wY \ot \wX
&
&\wY \ot \wY
&
&q^{-2}\wH \ot \wY +q^2\wY \ot \wH
\end{align*}

\noindent
\textbf{Braiding:} $(1-\sigma )(q^2+\sigma )=0$\\
$\ker (q^2+\sigma):$
$q\wH \ot \wX -q^{-1}\wX \ot \wH $,
$q^{-1}\wH \ot \wY -q\wY \ot \wH $,\\
\phantom{$\ker (q^2+\sigma):$}
$\wX \ot \wY -\wY \ot \wX $.

\newpage

\noindent
\textbf{9.\ Quantum tangent space {\boldmath $\X _\Ga $}:}
\begin{align*}
H&:=\frac{2}{q^{-2}-q^2}( K^4-1), & X&:=q^{-1/2}FK, & Y&:=q^{-1/2}KE
\end{align*}

\noindent
\textbf{Real forms:}
$\mc{O}(\mathrm{SU}_q(2))$,
$\mc{O}(\mathrm{SU}_q(1,1))$,
$\mc{O}(\mathrm{SL}_q(2,\mathbb{R}))$

\noindent
\textbf{Fund. repr.:}
$\displaystyle H=\frac{2}{q+q^{-1}}
\begin{pmatrix}q^{-1} & 0\\0 & -q\end{pmatrix}$,
$X=\begin{pmatrix}0 & 1\\0 & 0\end{pmatrix}$,
$Y=\begin{pmatrix}0 & 0\\1 & 0\end{pmatrix}$

\noindent
\textbf{Relations:}
\begin{align*}
q^2HX-q^{-2}XH&=2X\\
q^{-2}HY-q^2YH&=-2Y\\
q^{-1}XY-qYX&=\frac{q+q^{-1}}{2}H
\end{align*}

\noindent
\textbf{Module structure, differentials:}
\begin{align*}
(f^i_j)=
\begin{pmatrix}
K^4 & 0 & 0\\
0 & K^2 & 0\\
0 & 0 & K^2
\end{pmatrix}
\quad &
\begin{array}{rl}
\dif \wH =& \frac{-1-q^2}{2}\wX \wedge \wY \\
\dif \wX =& -2q^{-2}\wH \wedge \wX \\
\dif \wY =& 2q^2\wH \wedge \wY
\end{array}
\end{align*}

\noindent
\textbf{Pairing:}
\[
\begin{array}{r|ccccccccc}
 & (u^1_1)^2 & u^1_1u^1_2 & u^1_1u^2_1 & (u^1_2)^2 & u^1_2u^2_1 & u^1_2u^2_2
& (u^2_1)^2 & u^2_1u^2_2 & (u^2_2)^2\\
\hline
H & 2q^{-2} & 0 & 0 & 0 & 0 & 0 & 0 & 0 & -2q^2\\
X & 0 & 1 & 0 & 0 & 0 & q & 0 & 0 & 0\\
Y & 0 & 0 & 1 & 0 & 0 & 0 & 0 & q & 0
\end{array}
\]

\noindent
\textbf{Right ideal} {\boldmath $\mc{R}_\Ga $}\textbf{:}
$u^1_1+q^{-2}u^2_2-(1+q^{-2})$, $(u^1_2)^2$, $(u^2_1)^2$,
$u^1_2u^2_1$, $(u^1_1-1)u^1_2$, $(u^1_1-1)u^2_1$

\noindent
\textbf{Left-invariant symmetric 2-forms:}
\begin{align*}
&\wH \ot \wH
&
&\wX \ot \wX
&
&q^{-2}\wH \ot \wX +q^2\wX \ot \wH
\\
&q\wX \ot \wY +q^{-1}\wY \ot \wX
&
&\wY \ot \wY
&
&q^2\wH \ot \wY +q^{-2}\wY \ot \wH
\end{align*}

\noindent
\textbf{Braiding:} $(1-\sigma )(q^2+\sigma )=0$\\
$\ker (q^2+\sigma):$
$q^{-3}\wH \ot \wX -q^3\wX \ot \wH $,
$q^3\wH \ot \wY -q^{-3}\wY \ot \wH $,\\
\phantom{$\ker (q^2+\sigma):$}
$q^2\wX \ot \wY -q^{-2}\wY \ot \wX $.

\newpage

\noindent
\textbf{10.\ Quantum tangent space {\boldmath $\X _\Ga $}:}
\begin{align*}
H&:=\frac{2}{q^{-1}-q}( K^2-1), & X&:=q^{-1/2}FK, & Y&:=q^{-1/2}KE
\end{align*}

\noindent
\textbf{Real forms:}
$\mc{O}(\mathrm{SU}_q(2))$,
$\mc{O}(\mathrm{SU}_q(1,1))$,
$\mc{O}(\mathrm{SL}_q(2,\mathbb{R}))$

\noindent
\textbf{Fund. repr.:}
$\displaystyle H=\frac{2}{q+1}
\begin{pmatrix}1 & 0\\0 & -q\end{pmatrix}$,
$X=\begin{pmatrix}0 & 1\\0 & 0\end{pmatrix}$,
$Y=\begin{pmatrix}0 & 0\\1 & 0\end{pmatrix}$

\noindent
\textbf{Relations:}
\begin{align*}
qHX-q^{-1}XH&=2X\\
q^{-1}HY-qYH&=-2Y\\
q^{-1}XY-qYX+\frac{q-q^{-1}}{4}H^2&=H
\end{align*}

\noindent
\textbf{Module structure, differentials:}
\begin{align*}
(f^i_j)=
\begin{pmatrix}
K^2 & 0 & 0\\
0 & K^2 & 0\\
0 & 0 & K^2
\end{pmatrix}
\quad &
\begin{array}{rl}
\dif \wH =& -q\wX \wedge \wY \\
\dif \wX =& -2q^{-1}\wH \wedge \wX \\
\dif \wY =& 2q\wH \wedge \wY
\end{array}
\end{align*}

\noindent
\textbf{Pairing:}
\[
\begin{array}{r|ccccccccc}
 & (u^1_1)^2 & u^1_1u^1_2 & u^1_1u^2_1 & (u^1_2)^2 & u^1_2u^2_1 & u^1_2u^2_2
& (u^2_1)^2 & u^2_1u^2_2 & (u^2_2)^2\\
\hline
H & 2q^{-1} & 0 & 0 & 0 & 0 & 0 & 0 & 0 & -2q\\
X & 0 & 1 & 0 & 0 & 0 & q & 0 & 0 & 0\\
Y & 0 & 0 & 1 & 0 & 0 & 0 & 0 & q & 0
\end{array}
\]

\noindent
\textbf{Right ideal} {\boldmath $\mc{R}_\Ga $}\textbf{:}
$u^1_1+q^{-1}u^2_2-(1+q^{-1})$, $(u^1_2)^2$, $(u^2_1)^2$,
$u^1_2u^2_1$, $(u^1_1-1)u^1_2$, $(u^1_1-1)u^2_1$

\noindent
\textbf{Left-invariant symmetric 2-forms:}
\begin{align*}
&\wH \ot \wH +\frac{1-q^2}{4}\wX \ot \wY
&
&\wX \ot \wX
&
&q^{-1}\wH \ot \wX +q\wX \ot \wH
\\
&q\wX \ot \wY +q^{-1}\wY \ot \wX
&
&\wY \ot \wY
&
&q\wH \ot \wY +q^{-1}\wY \ot \wH
\end{align*}

\noindent
\textbf{Braiding:} $(1-\sigma )(q^2+\sigma )=0$\\
$\ker (q^2+\sigma):$
$\wH \ot \wX -\wX \ot \wH $,
$\wH \ot \wY -\wY \ot \wH $,
$\wX \ot \wY -\wY \ot \wX $.

\newpage

\noindent
\textbf{11.\ Quantum tangent space {\boldmath $\X _\Ga $}:}
\begin{align*}
H&:=\frac{2}{q-q^{-1}}(\coun _-K^2-1), & X&:=q^{-1/2}FK, & Y&:=q^{-1/2}KE
\end{align*}

\noindent
\textbf{Real forms:}
$\mc{O}(\mathrm{SU}_q(2))$,
$\mc{O}(\mathrm{SU}_q(1,1))$,
$\mc{O}(\mathrm{SL}_q(2,\mathbb{R}))$

\noindent
\textbf{Fund. repr.:}
$\displaystyle H=\frac{2}{1-q}
\begin{pmatrix}1 & 0\\0 & q\end{pmatrix}$,
$X=\begin{pmatrix}0 & 1\\0 & 0\end{pmatrix}$,
$Y=\begin{pmatrix}0 & 0\\1 & 0\end{pmatrix}$

\noindent
\textbf{Relations:}
\begin{align*}
qHX-q^{-1}XH&=-2X\\
q^{-1}HY-qYH&=2Y\\
q^{-1}XY-qYX+\frac{q-q^{-1}}{4}H^2&=-H
\end{align*}

\noindent
\textbf{Module structure, differentials:}
\begin{align*}
(f^i_j)=
\begin{pmatrix}
\coun _-K^2 & 0 & 0\\
0 & K^2 & 0\\
0 & 0 & K^2
\end{pmatrix}
\quad &
\begin{array}{rl}
\dif \wH =& q\wX \wedge \wY \\
\dif \wX =& 2q^{-1}\wH \wedge \wX \\
\dif \wY =& -2q\wH \wedge \wY
\end{array}
\end{align*}

\noindent
\textbf{Pairing:}
\[
\begin{array}{r|ccccccccc}
 & (u^1_1)^2 & u^1_1u^1_2 & u^1_1u^2_1 & (u^1_2)^2 & u^1_2u^2_1 & u^1_2u^2_2
& (u^2_1)^2 & u^2_1u^2_2 & (u^2_2)^2\\
\hline
H & -2q^{-1} & 0 & 0 & 0 & 0 & 0 & 0 & 0 & 2q\\
X & 0 & 1 & 0 & 0 & 0 & q & 0 & 0 & 0\\
Y & 0 & 0 & 1 & 0 & 0 & 0 & 0 & q & 0
\end{array}
\]

\noindent
\textbf{Right ideal} {\boldmath $\mc{R}_\Ga $}\textbf{:}
$u^1_1-q^{-1}u^2_2-(1-q^{-1})$, $(u^1_2)^2$, $(u^2_1)^2$,
$u^1_2u^2_1$, $(u^1_1-1)u^1_2$, $(u^1_1-1)u^2_1$

\noindent
\textbf{Left-invariant symmetric 2-forms:}
\begin{align*}
&\wH \ot \wH +\frac{1-q^2}{4}\wX \ot \wY
&
&\wX \ot \wX
&
&q^{-1}\wH \ot \wX +q\wX \ot \wH
\\
&q\wX \ot \wY +q^{-1}\wY \ot \wX
&
&\wY \ot \wY
&
&q\wH \ot \wY +q^{-1}\wY \ot \wH
\end{align*}

\noindent
\textbf{Braiding:} $(1-\sigma )(q^2+\sigma )=0$\\
$\ker (q^2+\sigma):$
$\wH \ot \wX -\wX \ot \wH $,
$\wH \ot \wY -\wY \ot \wH $,
$\wX \ot \wY -\wY \ot \wX $.

\newpage

\noindent
\textbf{Acknowledgement:}
The author would like to thank A.~Sch\"uler for stimulating discussions.


\end{document}